\documentclass[preprint,12pt]{elsarticle}
\usepackage{amssymb}
\usepackage[utf8]{inputenc}
\usepackage[T1]{fontenc}
\usepackage{amsfonts}
\usepackage{amssymb}
\usepackage{amsthm}
\usepackage{amsmath}
\usepackage{graphicx}
\usepackage{mathrsfs}

\newtheorem{theorem}{Theorem}

\newtheorem{definition}[theorem]{Definition}
\newtheorem{example}[theorem]{Example}

\newtheorem{lemma}[theorem]{Lemma}

\newtheorem{remark}[theorem]{Remark}

\journal{JDE}

\begin{document}

\begin{frontmatter}



\title{Singular Optimal Controls of Stochastic Recursive Systems and
Hamilton-Jacobi-Bellman Inequality}


\author{Liangquan Zhang\fnref{1}}

\address{{\small 1. School of Science, }\\
{\small Beijing University of Posts and Telecommunications, }\\
{\small Beijing 100876, China}}
\fntext[myfootnote]{The author acknowledges the financial support partly by the National Nature
Science Foundation of China (Grant No. 11701040, 61571059, 11871010 \&
61603049) and Innovation Foundation of BUPT for Youth (No. 500417024 \&
505018304). E-mail: xiaoquan51011@163.com}
\begin{abstract}
In this paper, we study the optimal singular controls for stochastic
recursive systems, in which the control has two components: the regular
control, and the singular control. Under certain assumptions, we establish
the dynamic programming principle for this kind of optimal singular controls
problem, and prove that the value function is a unique viscosity solution of
the corresponding Hamilton-Jacobi-Bellman inequality, in a given class of
bounded and continuous functions. At last, an example is given for
illustration.

\end{abstract}


\end{frontmatter}
\newpage

\section{Introduction}

\label{sect:1}Singular stochastic control problem is one of the most
important topics in fields of stochastic control. This problem was first
introduced by Bather and Chernoff \cite{BC} in 1967 by considering a
simplified model for the control of a spaceship. It was noted for this
special model that there was a connection between the singular control
problem and optimal stopping problem. This link was established through the
derivative of the value function of this initial singular control problem
and the value function of the corresponding optimal stopping problem. From
then on, it was extended by Ben\v{e}s, Shepp, Witzsenhausen (see \cite{BSW})
and Karatzas and Shreve (see \cite{K1, K2, KS1, KS2, KS3}). In contrast to
classical control problems, singular control problems admit both of the
continuity of the cumulative displacement of the state caused by control and
the jump of one in impulsive control problems, between which it is either
constant or absolutely continuous.

Recall the classical singular control problems: The state process is driven
by the following $n$-dimensional stochastic differential equation%
\begin{equation}
\left\{
\begin{array}{lll}
\mathrm{d}X_{s}^{t,x;v,\xi } & = & b\left( s,X_{s}^{t,x;v,\xi },v_{s}\right)
\mathrm{d}s+\sigma \left( s,X_{s}^{t,x;v,\xi },v_{s}\right) \mathrm{d}%
W_{s}+G_{s}\mathrm{d}\xi _{s}, \\
X_{t}^{t,x;v,\xi } & = & x,\qquad 0\leq t\leq s\leq T,%
\end{array}%
\right.  \label{SDE1}
\end{equation}%
on some filtered probability space $\left( \Omega ,\mathscr{F},P\right) $,
where $b\left( \cdot ,\cdot ,\cdot \right) :\left[ 0,T\right] \times \mathbb{%
R}^{n}\times \mathbb{R}^{k}\rightarrow \mathbb{R}^{n},$ $\sigma \left( \cdot
,\cdot ,\cdot \right) :\left[ 0,T\right] \times \mathbb{R}^{n}\times \mathbb{%
R}^{k}\rightarrow \mathbb{R}^{n\times d},$ $G\left( \cdot \right) :\left[ 0,T%
\right] \rightarrow \mathbb{R}^{n\times m}$ are given deterministic
functions, $\left( W_{s}\right) _{s\geq 0}$ is an $d$-dimensional Brownian
motion, $\left( x,t\right) $ are initial time and state, $v\left( \cdot
\right) :\left[ 0,T\right] \rightarrow \mathbb{R}^{k}$ is classical control
process, and $\xi \left( \cdot \right) :\left[ 0,T\right] \rightarrow
\mathbb{R}^{m}$, with nondecreasing left-continuous with right limits stands
for the singular control.

The aim is to minimize the cost functional:%
\begin{equation}
J\left( t,x;v,\xi \right) =\mathbb{E}\left[ \int_{t}^{T}l\left(
s,X_{s}^{t,x;v,\xi },v_{s}\right) \mathrm{d}s+\int_{t}^{T}K_{s}\mathrm{d}\xi
_{s}\right] ,  \label{cost1}
\end{equation}%
where
\begin{eqnarray*}
l\left( \cdot ,\cdot ,\cdot \right) &:&\left[ 0,T\right] \times \mathbb{R}%
^{n}\times \mathbb{R}^{k}\rightarrow \mathbb{R}, \\
K\left( \cdot \right) &:&\left[ 0,T\right] \rightarrow \mathbb{R}%
_{+}^{m}\triangleq \left\{ x\in \mathbb{R}^{m}:x_{i}\geq 0,\text{ }%
i=1,\ldots m\right\}
\end{eqnarray*}
are given deterministic functions, where $l\left( \cdot \right) $ represents
the running cost tare of the problem and $K\left( \cdot \right) $ the cost
rate of applying the singular control.

There are four approaches employed: the first one is based on the theory of
partial differential equations and on variational arguments, and can be
found in the works of Alvarez \cite{A1, A2}, Chow, Menaldi, and Robin \cite%
{CMR}, Karatzas \cite{K2}, Karatzas and Shreve \cite{KS3}, and Menaldi and
Taksar \cite{MT}. The second one is related to probabilistic methods; see
Baldursson \cite{B1}, Boetius \cite{Bo1, Bo2}, Boetius and Kohlmann \cite%
{BK1}, El Karoui and Karatzas \cite{EK1, EK2}, Karatzas \cite{K1}, and
Karatzas and Shreve \cite{KS1, KS2}. Third, the DPP, has been studied in a
general context, for example, by Boetius \cite{Bo2}, Haussmann and Suo \cite%
{HS2}, Fleming and Soner \cite{FS}, and Zhu \cite{Zhu}. At last the maximum
principle for optimal singular controls (see, for example, Cadenillas and
Haussmann \cite{CH}, Dufour and Miller \cite{DM2}, Dahl and \O ksendal \cite%
{DO} see references therein). The existence for optimal singular control can
be found in Haussmann and Suo \cite{HS2} and Dufour and Miller \cite{DM1}
via different approaches.

It is necessary to point out that singular control problems are largely used
in diverse areas such as mathematical finance (see Baldursson and Karatzas
\cite{BK}, Chiarolla and Haussmann \cite{CHF}, Kobila \cite{KOB}, and
Karatzas and Wang \cite{KW}, Davis and Norman \cite{DN}), manufacturing
systems (see, Shreve, Lehoczky, and Gaver \cite{SLG}), and queuing systems
(see Martins and Kushner \cite{MK}). Particularly, the application of H-J-B
inequality in finance can be seen in Pag\`{e}s and Possama\"{\i} \cite{PPs},
which is employed to investigate the bank monitoring incentives.

As is well known for the classical stochastic control problems, DPP is
satisfied. Moreover, if the value function has appropriate regularity, it
admits a second-order nonlinear partial differential equation (H-J-B
equation) (see Fleming and Rishel \cite{FR} and Lions \cite{Lion} reference
therein). In the frame work of singular stochastic control, the H-J-B
equation becomes a second-order variational inequality (see Fleming and
Soner \cite{FS}, Haussmannand, Suo \cite{HS1}).

Wang \cite{W} firstly introduces and studies a class of singular control
problems with recursive utility, where the cost function is determined by a
backward stochastic differential equation (BSDE in short). More preciously,
the cost functional is defined by%
\begin{equation}
J\left( t,x;\xi \right) =Y_{t}^{t,x;\xi },  \label{cost2}
\end{equation}%
where $Y_{t}^{t,x;v,\xi }$ is determined by
\begin{equation}
\left\{
\begin{array}{lll}
\mathrm{d}X_{s}^{t,x;\xi } & = & \left( aX_{s}^{t,x;\xi }+b\right) \mathrm{d}%
s+\sigma \mathrm{d}W_{s}+\mathrm{d}\xi _{s}, \\
\mathrm{d}Y_{s}^{t,x;\xi } & = & -f\left( Y_{s}^{t,x;\xi },Z_{s}^{t,x;\xi
},v_{s}\right) -\mu Z_{s}^{t,x;\xi }\mathrm{d}s+Z_{s}^{t,x;\xi }\mathrm{d}%
W_{s}-\mathrm{d}\xi _{s}, \\
X_{t}^{t,x;\xi } & = & x,\text{ }Y_{T}^{t,x;\xi }=0,\qquad 0\leq t\leq s\leq
T.%
\end{array}%
\right.  \label{Wang}
\end{equation}%
The value function is defined by
\begin{equation*}
u\left( t,x\right) =\inf_{\xi \in \mathcal{V}}J\left( t,x;\xi \right) ,
\end{equation*}%
where $\mathcal{V}$ denotes the set of admissible controls $\xi \left( \cdot
\right) $ which are $\mathscr{F}_{t}$-adapted processes from $\left[ 0,T%
\right] $ into $\mathbb{R}$, left continuous having right limits,
nonnegative and nondecreasing, moreover $\xi _{t}=0$. Under certain
assumptions, the author proved that the value function is a nonnegative,
convex solution of the H-J-B equation:%
\begin{equation*}
\min \left\{ \frac{\partial }{\partial t}u+\frac{1}{2}\sigma ^{2}\frac{%
\partial ^{2}}{\partial x^{2}}u+\left( ax+b\right) \frac{\partial }{\partial
x}u,u_{x}+1\right\} =0,\text{ }u\left( T,x\right) =0.
\end{equation*}%
As a matter of fact, the class of stochastic control problems with recursive
utility was firstly considered by Peng \cite{P2} where he considered an
absolutely continuous stochastic control problem with recursive utility, in
which the cost function is defined by a BSDE. As is well known, BSDE has
been widely used in mathematical finance and differential games, partial
differential equation, \textit{etc} (see \cite{KPQ, P1, P2, HL, WY, LP}).
Since then, there exists a huge literature to focus on the stochastic
recursive optimal control problems (see reference therein).

Motivated by above work, in this paper, we consider singular controls
problem of a general framework like:%
\begin{equation}
\left\{
\begin{array}{lll}
\mathrm{d}X_{s}^{t,x;v,\xi } & = & b\left( s,X_{s}^{t,x;v,\xi },v_{s}\right)
\mathrm{d}s+\sigma \left( s,X_{s}^{t,x;v,\xi },v_{s}\right) \mathrm{d}W_{s}+G%
\mathrm{d}\xi _{s}, \\
\mathrm{d}Y_{s}^{t,x;v,\xi } & = & -f\left( t,X_{s}^{t,x;v,\xi
},Y_{s}^{t,x;v,\xi },Z_{s}^{t,x;v,\xi },v_{s}\right) \mathrm{d}%
s+Z_{s}^{t,x;v,\xi }\mathrm{d}W_{s}-K\mathrm{d}\xi _{s}, \\
X_{t}^{t,x;v,\xi } & = & x,\text{ }Y_{T}^{t,x;v,\xi }=\Phi \left(
X_{T}^{t,x;v,\xi }\right) ,\qquad 0\leq t\leq s\leq T,%
\end{array}%
\right.  \label{FBSDE1}
\end{equation}%
with the similar cost functional%
\begin{equation}
J\left( t,x;\xi \right) =Y_{t}^{t,x;v,\xi }.
\end{equation}%
Note that $G$ and $K$ are deterministic matrices in Eq. (\ref{FBSDE1}). The
justification will be given in Remark \ref{GK} below. Using some properties
of the BSDE and analysis technique, we expand the extension of the dynamic
programming principle of the recursive control problem in \cite{BBP1, CP,
CW, LW, P2, SXZ, W, ZZ} to the singular controls case. And then, we show
that, provided the problem is formulated within a Markovian framework, the
value function is a unique viscosity solution of the problem for one kind of
nonlinear H-J-B inequality, in a given class of bounded and continuous
functions. Some characteristics of the value function of the control problem
are obtained as well. A concrete example is provided to validate our
theoretic results.

To the best of our knowledge, such singular optimal controls of FBSDEs (\ref%
{FBSDE1}) via DPP have not been studied before. On the one hand, our paper
generalizes the DPP of the forward controls problem by Haussmann and Suo
\cite{HS2} to recursive situation, backward one in Peng \cite{P2, P3}, Wu
and Yu \cite{WY} to the \textit{singular} control case, and admits the
following advantages and improvements. First, in Peng \cite{P2, P3}, Wu and
Yu \cite{WY}, the recursive cost function does not have the singular
control. In our paper, the optimization problem has singular control for the
cost function which is described by the solution of BSDE. In contrast to Wu
and Yu \cite{WY}, our \textit{variational inequality} involves the
derivative of the value function. As a matter of fact, Example 1.7 in the
classical literature by Crandall, Ishii and Lions \cite{CIL} displays such a
kind of obstacle and gradient constraint problems. Thus, the \textit{purpose}
of this paper is to response this item and to deepen and enrich this
research. Besides, through the value function, we can define the so called
\textit{inaction region}, in which the status process is continuous and the
optimal singular control should not produce any jump. Our ideas and
methodology of proof in our paper for DPP come from Haussmannand, Suo \cite%
{HS2}, Peng \cite{P2, P3}, Wu and Yu \cite{WY}. Second, the method of proof
in our paper is mainly based on elementary mathematics analysis technique
and the properties of BSDE. Wang \cite{W} studies a particular situation.
However, the FBSDEs in Wang \cite{W} do not contain the regular control and
the generator is not general case. The solution there belongs to Sobolev
space. Nonetheless, in classical stochastic control theory, viscosity
solution is usually employed to study H-J-B equation. In this paper, we
adopt the approach developed by Fleming and Soner \cite{FS} to get the
uniqueness proof for the viscosity solution of corresponding H-J-B
inequality.

The rest of this paper is organized as follows: After some preliminaries in
Section \ref{sect:2}, we are devoted the third section to the DPP for the
optimal singular controls problem and prove that the value function of the
singular controls problem is the unique viscosity solution of the
corresponding H-J-B inequality. We conclude with a concrete example in
Section \ref{sect:4}. Finally, some proofs of lemmas are scheduled in the
Appendix.

\section{Formulation of the Problem}

\label{sect:2}

Throughout this paper, we denote by $\mathbb{R}^{n}$ the space of $n$%
-dimensional Euclidean space, by $\mathbb{R}^{n\times d}$ the space the
matrices with order $n\times d$. The matrix transpose is denoted by $\top .$
Let $(\Omega ,\mathscr{F},\{\mathscr{F}_{t}\}_{t\geq 0},P)$ be a complete
filtered probability space on which a $d$-dimensional standard Brownian
motion $W(\cdot )$ is defined, with $\{\mathscr{F}_{t}\}_{t\geq 0}$ being
its natural filtration, augmented by all the $P$-null sets. Given a subset $%
U $ (compact) of $\mathbb{R}^{k},$ we will denote $\mathcal{U}$ the class of
measurable, adapted processes $\left( v,\xi \right) :\left[ 0,T\right]
\times \Omega \rightarrow U\times \left[ 0,\infty \right) ^{m},$ with $\xi $
nondecreasing left-continuous with right limits and $\xi _{0}=0$, moreover, $%
\mathbb{E}\left[ \sup\limits_{0\leq t\leq T}\left\vert v_{t}\right\vert
^{2}+\left\vert \xi _{T}\right\vert ^{2}\right] <+\infty .$ For each $t>0$,
we denote by $\left\{ \mathscr{F}_{s}^{t},t\leq s\leq T\right\} $ the
natural filtration of the Brownian motion $\left\{ W_{s}-W_{t},\text{ }t\leq
s\leq T\right\} $, augmented by the $P$-null sets of $\mathscr{F}$. In what
follows, $C$ represents a generic constant, which can be different from line
to line.

We now introduce the following spaces of processes:
\begin{align*}
\mathcal{S}^{2}(0,T;\mathbb{R})\triangleq & \left\{ \mathbb{R}^{n}\text{%
-valued }\mathscr{F}_{t}\text{-adapted process }\phi (t)\text{; }\mathbb{E}%
\left[ \sup\limits_{0\leq t\leq T}\left\vert \phi _{t}\right\vert ^{2}\right]
<\infty \right\} , \\
\mathcal{M}^{2}(0,T;\mathbb{R})\triangleq & \left\{ \mathbb{R}^{n}\text{%
-valued }\mathscr{F}_{t}\text{-adapted process }\varphi (t)\text{; }\mathbb{E%
}\left[ \int_{0}^{T}\left\vert \varphi _{t}\right\vert ^{2}\mbox{\rm d}t%
\right] <\infty \right\} ,
\end{align*}%
and denote $\mathcal{N}^{2}\left[ 0,T\right] =\mathcal{S}^{2}(0,T;\mathbb{R}%
^{n})\times \mathcal{S}^{2}(0,T;\mathbb{R})\times \mathcal{M}^{2}(0,T;%
\mathbb{R}^{n}).$ Clearly, $\mathcal{N}^{2}\left[ 0,T\right] $ forms a
Banach space.

For any $v\left( \cdot \right) \times \xi \left( \cdot \right) \in \mathcal{U%
},$ we study the stochastic control systems governed by FBSDE of the
following type with two adapted control processes:
\begin{equation}
\left\{
\begin{array}{rcl}
\mathrm{d}X_{s}^{t,x;v,\xi } & = & b\left( s,X_{s}^{t,x;v,\xi },v_{s}\right)
\mathrm{d}s+\sigma \left( s,X_{s}^{t,x;v,\xi },v_{s}\right) \mathrm{d}W_{s}+G%
\mathrm{d}\xi _{s}, \\
\mathrm{d}Y_{s}^{t,x;v,\xi } & = & -f\left( s,X_{s}^{t,x;v,\xi
},Y_{s}^{t,x;v,\xi },Z_{s}^{t,x;v,\xi },v_{s}\right) \mathrm{d}%
s+Z_{s}^{t,x;v,\xi }\mathrm{d}W_{s}-K\mathrm{d}\xi _{s}, \\
X_{t}^{t,x;v,\xi } & = & x,\text{ }Y_{T}^{t,x;v,\xi }=\Phi \left(
X_{T}^{t,x;v,\xi }\right) ,\qquad 0\leq t\leq s\leq T,%
\end{array}%
\right.  \label{eq1}
\end{equation}%
\newline
where
\begin{eqnarray*}
b &:&\left[ 0,T\right] \times \mathbb{R}^{n}\times \mathbb{R}^{k}\rightarrow
\mathbb{R}^{n}, \\
\sigma &:&\left[ 0,T\right] \times \mathbb{R}^{n}\times \mathbb{R}%
^{k}\rightarrow \mathbb{R}^{n\times d}, \\
f &:&\left[ 0,T\right] \times \mathbb{R}^{n}\times \mathbb{R}\times \mathbb{R%
}^{d}\times \mathbb{R}^{k}\rightarrow \mathbb{R}, \\
\Phi &:&\mathbb{R}^{n}\rightarrow \mathbb{R}.
\end{eqnarray*}

\begin{definition}
For any $v\left( \cdot \right) \times \xi \left( \cdot \right) \in \mathcal{U%
}$, a triple of processes
\begin{equation*}
\left( X_{\cdot }^{t,x;v,\xi },Y_{\cdot }^{t,x;v,\xi },Z_{\cdot }^{t,x;v,\xi
}\right) :\left[ 0,T\right] \times \Omega \rightarrow \mathbb{R}^{n}\times
\mathbb{R}\times \mathbb{R}^{d}
\end{equation*}%
is called an adapted solution of the FBSDEs \emph{(\ref{eq1})}, if $\left(
X_{\cdot }^{t,x;v,\xi },Y_{\cdot }^{t,x;v,\xi },Z_{\cdot }^{t,x;v,\xi
}\right) \in \mathcal{N}^{2}\left[ 0,T\right] ,$ and it satisfies \emph{(\ref%
{eq1})}, $P$-almost surely.
\end{definition}

For simplicity, we adopt the following assumptions taken from Peng \cite{P3}.

\begin{description}
\item[(A1)] $b$, $\sigma $ are uniformly continuous in $\left( s,x,u\right) $%
, and there exits a positive constant $C$ such that for all $s\in \left[ 0,T%
\right] $, $x^{1},$ $x^{2}\in \mathbb{R}^{n}$, $u^{1},$ $u^{2}\in \mathbb{R}%
^{k},$%
\begin{eqnarray*}
&&\left\vert b\left( s,x^{1},u^{1}\right) -b\left( s,x^{2},u^{2}\right)
\right\vert +\left\vert \sigma \left( s,x^{1},u^{1}\right) -\sigma \left(
s,x^{2},u^{2}\right) \right\vert \\
&\leq &C\left( \left\vert x^{1}-x^{2}\right\vert +\left\vert
u^{1}-u^{2}\right\vert \right) ,
\end{eqnarray*}%
and%
\begin{equation*}
\left\vert b\left( s,x,u\right) \right\vert +\left\vert \sigma \left(
s,x,u\right) \right\vert \leq C\left( 1+\left\vert x\right\vert \right) .
\end{equation*}

\item[(A2)] $f,$ $\Phi $ are uniformly continuous in $\left( s,x,y,u\right) $
and there exists a constant $C>0$ such that for all $s\in \left[ 0,T\right] $%
, $x^{1},$ $x^{2}\in \mathbb{R}^{n}$, $y^{1},$ $y^{2}\in \mathbb{R},$ $%
z^{1}, $ $z^{2}\in \mathbb{R}^{d},$ $u^{1},$ $u^{2}\in \mathbb{R}^{k},$%
\begin{eqnarray*}
&&\left\vert f\left( s,x^{1},y^{1},z^{1},u^{1}\right) -f\left(
s,x^{2},y^{2},z^{2},u^{2}\right) \right\vert +\left\vert \Phi \left(
x^{1}\right) -\Phi \left( x^{2}\right) \right\vert \\
&\leq &C\left( \left\vert x^{1}-x^{2}\right\vert +\left\vert
y^{1}-y^{2}\right\vert +\left\vert z^{1}-z^{2}\right\vert +\left\vert
u^{1}-u^{2}\right\vert \right) ,
\end{eqnarray*}%
moreover,%
\begin{equation*}
\left\vert f\left( s,x,0,0,u\right) \right\vert +\left\vert \Phi \left(
x\right) \right\vert \leq C\left( 1+\left\vert x\right\vert \right) .
\end{equation*}

\item[(A3)] $G_{n\times m}$ and $K_{1\times m}$ are deterministic matrices.
Moreover, postulate that there exists $k_{0}$ such that $K^{i}>k_{0}>0,$ $%
1\leq i\leq m,$ where $K^{i}$ denotes the $i$th coordinate of $K.$
\end{description}

\begin{remark}
\label{GK} In this paper, we assume that $G_{n\times m}$ and $K_{1\times m}$
are deterministic matrices. On the one hand, from the derivations in Theorem
5.1 of \cite{HS2}, also in Lemma \ref{Property} below, in our paper, it is
convenient to show the \textquotedblleft inaction\textquotedblright\ region
for singular control; On the other hand, we may regard $Y_{s}^{t,x;v,\xi
}+G\xi _{s}$ together as a solution, in this way, we are able to apply the
classical It\^{o}'s formula, avoiding the appearance of jump. We believe
these assumptions can be removed properly, but at present, we consider
constant only in our paper for our aim. Whilst in order to get the
uniqueness of the solution to H-J-B inequality (\ref{pde1}), we add the
assumption $K^{i}>k_{0}>0,$ $1\leq i\leq m.$ More details, see Theorem \ref%
{t4} below.
\end{remark}

Under above assumptions (A1)-(A3), for any $v\left( \cdot \right) \times
\xi \left( \cdot \right) \in \mathcal{U}$, it is easy to check that FBSDEs (%
\ref{eq1}) admit a unique $\mathcal{F}_{t}$-adapted solution denoted by the
triple $(X_{\cdot }^{t,x;v,\xi },Y_{\cdot }^{t,x;v,\xi },Z_{\cdot
}^{t,x;v,\xi })\in \mathcal{N}^{2}\left[ 0,T\right] $ (See Pardoux and Peng
\cite{PP1}).

Like Peng \cite{P3}, given any control processes $v\left( \cdot \right)
\times \xi \left( \cdot \right) \in \mathcal{U}$, we introduce the following
cost functional:
\begin{equation}
J(t,x;v\left( \cdot \right) ,\xi \left( \cdot \right) )=\left.
Y_{s}^{t,x;v,\xi }\right\vert _{s=t},\qquad \left( t,x\right) \in \left[ 0,T%
\right] \times \mathbb{R}^{n}.  \label{c1}
\end{equation}%
We are interested in the \emph{value function} of the stochastic optimal
control problem%
\begin{eqnarray}
u\left( t,x\right) &=&J(t,x;\hat{v}\left( \cdot \right) ,\hat{\xi}\left(
\cdot \right) )  \notag \\
&=&ess\inf_{v\left( \cdot \right) \times \xi \left( \cdot \right) \in
\mathcal{U}}J\left( t,x;v\left( \cdot \right) ,\xi \left( \cdot \right)
\right) ,\text{ }\left( t,x\right) \in \left[ 0,T\right] \times \mathbb{R}%
^{n}.  \label{value}
\end{eqnarray}%
Since the value function (\ref{value}) is defined by the solution of
controlled BSDE (\ref{eq1}), so from the existence and uniqueness, $u$ is
\emph{well-defined}.

Any $\hat{v}\left( \cdot \right) \times \hat{\xi}\left( \cdot \right) \in
\mathcal{U}$ satisfying (\ref{value}) is called an optimal control pair of
optimal singular problem, and the corresponding state processes, denoted by $%
\left( X_{\cdot }^{t,x;\hat{v},\hat{\xi}},Y_{\cdot }^{t,x;\hat{v},\hat{\xi}%
},Z_{\cdot }^{t,x;\hat{v},\hat{\xi}}\right) $, is called optimal state
process. We also refer to $(X_{\cdot }^{t,x;\hat{v},\hat{\xi}},Y_{\cdot
}^{t,x;\hat{v},\hat{\xi}},Z_{\cdot }^{t,x;\hat{v},\hat{\xi}},\hat{v}\left(
\cdot \right) ,\hat{\xi}\left( \cdot \right) )$ as an optimal 5-tuple of
optimal singular problem$.$

\begin{definition}[Optimal Control]
Any admissible controls $\hat{u}(\cdot )\times \hat{\xi}\left( \cdot \right)
\in \mathcal{U},$ are called optimal$,$ if $\hat{u}(\cdot )\times \hat{\xi}%
\left( \cdot \right) $ attains the minimum of $J(u(\cdot )\times \xi \left(
\cdot \right) ).$
\end{definition}

Reader interested in stochastic optimal control problems can refer to Yong
and Zhou \cite{YZ}. We shall recall the following basic result on BSDE. We
begin with the well-known comparison theorem (see Barles, Buckdahn, and
Pardoux \cite{BBP}, Proposition 2.6).

\begin{lemma}[Comparison theorem]
\label{comp}Let $\left( y^{i},z^{i}\right) ,$ $i=1,2,$ be the solution to
the following
\begin{equation}
y^{i}\left( t\right) =\xi ^{i}+\int_{t}^{T}f^{i}\left(
s,y_{s}^{i},z_{s}^{i}\right) \mathrm{d}s-\int_{t}^{T}z_{s}^{i}\mathrm{d}%
W_{s},  \label{estbdsde}
\end{equation}%
where $\mathbb{E}\left[ \left\vert \xi ^{i}\right\vert ^{2}\right] <\infty ,$
$f^{i}\left( s,y^{i},z^{i}\right) $ satisfies the conditions \emph{(A2), }$%
i=1,2$\textsl{.} Under assumption \emph{(A2)}, BSDE \emph{(\ref{estbdsde}%
)} admits a unique adapted solution $\left( y^{i},z^{i}\right) ,$
respectively, for $i=1,2$. Furthermore, if

\noindent (i) $\xi ^{1}\geq \xi ^{2},$ a.s.;

\noindent (ii) $f^{1}\left( t,y,z\right) \geq f^{2}\left( t,y,z\right) ,$
a.e., for any $\left( t,y,z\right) \in \left[ 0,T\right] \times \mathbb{R}%
\times \mathbb{R}^{d}.$\newline
Then, we have: $y_{t}^{1}\geq y_{t}^{2}$, a.s., for all $t\in \left[ 0,T%
\right] .$
\end{lemma}

\section{Dynamic Programming Principle}

\label{sect:3}

In this section, we shall establish the DPP, which plays an important role
to derive a new H-J-B inequality. Furthermore, the notion of viscosity
solution to dynamic programming inequality will be introduced a moment
later. Before this, serval technique lemmas will be presented. From these
series of properties of value function (\ref{value}), we are able to prove
it is a viscosity solution.

\begin{definition}[Backward Semigroup]
\label{backward}For every $s\in \left[ t,t_{1}\right] $ with $t\leq
t_{1}\leq T$ and $\eta \in L^{2}\left( \Omega ,\mathscr{F}_{t_{1}},P;\mathbb{%
R}\right) ,$ we denote%
\begin{equation*}
\mathcal{G}_{s,t_{1}}^{t,x;v,\xi }\left[ \eta \right] =Y^{t,x;v,\xi }\left(
s\right) ,
\end{equation*}%
where%
\begin{equation*}
\left\{
\begin{array}{rcl}
\mathrm{d}X_{s}^{t,x;v,\xi } & = & b\left( s,X_{s}^{t,x;v,\xi },v_{s}\right)
\mathrm{d}s+\sigma \left( s,X_{s}^{t,x;v,\xi },v_{s}\right) \mathrm{d}W_{s}+G%
\mathrm{d}\xi _{s}, \\
\mathrm{d}Y_{s}^{t,x;v,\xi } & = & -f\left( t,X_{s}^{t,x;v,\xi
},Y_{s}^{t,x;v,\xi },Z_{s}^{t,x;v,\xi },v_{t}\right) \mathrm{d}%
s+Z_{s}^{t,x;v,\xi }\mathrm{d}W_{s}-K\mathrm{d}\xi _{s}, \\
X_{t}^{t,x;v,\xi } & = & x,\text{ }Y_{t_{1}}^{t,x;v,\xi }=\eta ,\qquad 0\leq
t\leq s\leq T.%
\end{array}%
\right.
\end{equation*}
\end{definition}

We start from following:

\begin{definition}
Given any $t\in \left[ 0,T\right] ,$ a sequence $\left\{ A_{i}\right\}
_{i=1}^{N}\subset \mathscr{F}_{t}$ is called a partition of $\left( \Omega ,%
\mathscr{F}_{t}\right) $ if $\cup _{i=1}^{N}A_{i}=\Omega $ and $A_{i}\cap
A_{j}=\phi $, whenever $i\neq j.$
\end{definition}

\begin{lemma}
\label{state3}Under assumptions \emph{(A1)-(A3)}, the value function $%
u\left( t,x\right) $ defined in \emph{(\ref{value})} is a deterministic
function.
\end{lemma}

\paragraph{Proof}

We adopt the idea from Peng \cite{P3}. In fact, there exists $\left(
v^{l},\xi ^{l}\right) _{l\geq 1}$ such that $u\left( t,x\right) =\inf_{l\geq
1}J\left( t,x,v^{l},\xi ^{l}\right) ,$ a.s. For any $\left( \left( v,\xi
\right) ,\left( v^{\prime },\xi ^{\prime }\right) \right) \in \mathcal{U}$,
we define
\begin{equation*}
\left( \left( v,\xi \right) \wedge \left( v^{\prime },\xi ^{\prime }\right)
\right) _{s}=\left\{
\begin{array}{l}
0,\text{ }s\in \left[ 0,t\right] ; \\
\left( v,\xi \right) _{s},\text{ }s\in \left( t,T\right] \text{ on }\left\{
J\left( t,x,v,\xi \right) \leq J\left( t,x,v^{\prime },\xi ^{\prime }\right)
\right\} ; \\
\left( v^{\prime },\xi ^{\prime }\right) _{s},\text{ }s\in \left( t,T\right]
\text{ on }\left\{ J\left( t,x,v,\xi \right) \geq J\left( t,x,v^{\prime
},\xi ^{\prime }\right) \right\} .%
\end{array}%
\right.
\end{equation*}%
Hence, $\left( \left( v,\xi \right) \vee \left( v^{\prime },\xi ^{\prime
}\right) \right) \in \mathcal{U}$, moreover, it follows that
\begin{equation*}
J\left( t,x,\left( v,\xi \right) \wedge \left( v^{\prime },\xi ^{\prime
}\right) \right) \leq J\left( t,x,v,\xi \right) \wedge J\left( t,x,v^{\prime
},\xi ^{\prime }\right) .
\end{equation*}%
Define $\left( \bar{v}^{1},\bar{\xi}^{1}\right) :=\left( v^{1},\xi
^{1}\right) ,$ $\left( \bar{v}^{l},\bar{\xi}^{l}\right) :=\left( \bar{v}%
^{l-1},\bar{\xi}^{l-1}\right) \wedge \left( v^{l},\xi ^{l}\right) .$ Then,
we have%
\begin{equation*}
u\left( t,x\right) =\lim_{l\uparrow \infty }\downarrow J\left( t,x,\bar{v}%
^{l},\bar{\xi}^{l}\right) .
\end{equation*}%
Suppose that
\begin{equation*}
\mathbb{E}\left[ \left\vert u\left( t,x\right) -J\left( t,x,\bar{v}^{l},\bar{%
\xi}^{l}\right) \right\vert ^{2}\right] \leq \frac{1}{l},\text{ }l\geq 1.
\end{equation*}%
For all $l\geq 1,$ there exists $\left( \tilde{v},\tilde{\xi}\right) $ where
\begin{eqnarray*}
\tilde{v}^{l} &=&\sum_{r,i,j}^{N_{l}-1}\bar{v}_{i,j,r}\mathbf{1}%
_{A_{j}^{l}\cap B_{ir}^{l}}\mathbf{1}_{\left( t_{i}^{l},t_{i+1}^{l}\right] },%
\text{ }t=t_{0}^{l}\leq \cdots t_{N_{l}}^{l}=T,\text{ }\bar{v}_{i,j,r}\in U,
\\
\tilde{\xi}^{l} &=&\sum_{r,i,j}^{N_{l}-1}\bar{\xi}_{i,j,r}\mathbf{1}%
_{A_{j}^{l}\cap B_{ir}^{l}}\mathbf{1}_{\left( t_{i}^{l},t_{i+1}^{l}\right] },%
\text{ }t=t_{0}^{l}\leq \cdots t_{N_{l}}^{l}=T,\text{ }\bar{\xi}_{i,j,r}\in
\left( \left[ 0,\infty \right) \right) ^{m},
\end{eqnarray*}%
where $\left\{ A_{j}^{l}\right\} _{0\leq j\leq N_{l}-1}$ is a partition of $%
\left( \Omega ,\mathscr{F}_{t}\right) ,$ $\left\{ B_{ir}^{l}\right\} _{0\leq
r\leq N_{l}-1}$ is a partition of $\left( \Omega ,\mathscr{F}%
_{t_{i}^{l}}^{t}\right) ,$ such that
\begin{equation*}
\mathbb{E}\left[ \int_{t}^{T}\left( \left\vert \bar{v}_{s}^{l}-\tilde{v}%
_{s}^{l}\right\vert ^{2}+\left\vert \bar{\xi}_{s}^{l}-\tilde{\xi}%
_{s}^{l}\right\vert ^{2}\right) \mathrm{d}s++\left\vert \bar{\xi}_{T}^{l}-%
\tilde{\xi}_{T}^{l}\right\vert ^{2}\right] \leq \frac{C}{l},\text{ }l\geq 1.
\end{equation*}%
From (\ref{EE4}), it derives that
\begin{eqnarray}
&&\mathbb{E}\left[ \left\vert J\left( t,x,\tilde{v}^{l},\tilde{\xi}%
^{l}\right) -J\left( t,x,\bar{v}^{l},\bar{\xi}^{l}\right) \right\vert ^{2}%
\right]  \notag \\
&\leq &C\mathbb{E}\left[ \int_{t}^{T}\left( \left\vert \bar{v}_{s}^{l}-%
\tilde{v}_{s}^{l}\right\vert ^{2}+\left\vert \bar{\xi}_{s}^{l}-\tilde{\xi}%
_{s}^{l}\right\vert ^{2}\right) \mathrm{d}s+\left\vert \bar{\xi}_{T}^{l}-%
\tilde{\xi}_{T}^{l}\right\vert ^{2}\right]  \notag \\
&\leq &\frac{C}{l},\text{ }l\geq 1.  \label{S1}
\end{eqnarray}%
From (\ref{S1}), we get%
\begin{equation}
\mathbb{E}\left[ \left\vert u\left( t,x\right) -J\left( t,x,\bar{v}^{l},\bar{%
\xi}^{l}\right) \right\vert ^{2}\right] \leq \frac{C}{l}.  \label{coner1}
\end{equation}%
Now we put%
\begin{eqnarray*}
\tilde{v}^{l,j} &=&\sum_{r,i,j}^{N_{l}-1}\bar{v}_{i,j,r}\mathbf{1}%
_{A_{j}^{l}\cap B_{ir}^{l}}\mathbf{1}_{\left( t_{i}^{l},t_{i+1}^{l}\right] },%
\text{ }t=t_{0}^{l}\leq \cdots t_{N_{l}}^{l}=T,\text{ }\bar{v}_{i,j,r}\in U,
\\
\tilde{\xi}^{l,j} &=&\sum_{r,i,j}^{N_{l}-1}\bar{\xi}_{i,j,r}\mathbf{1}%
_{A_{j}^{l}\cap B_{ir}^{l}}\mathbf{1}_{\left( t_{i}^{l},t_{i+1}^{l}\right] },%
\text{ }t=t_{0}^{l}\leq \cdots t_{N_{l}}^{l}=T,\text{ }\bar{\xi}_{i,j,r}\in
\left( \left[ 0,\infty \right) \right) ^{m}.
\end{eqnarray*}%
They are clearly $\mathscr{F}_{s}^{t}$-adapted $U\times \left( \left[
0,\infty \right) \right) ^{m}$ processes. Now consider the following%
\begin{equation*}
\left\{
\begin{array}{lll}
\mathbf{1}_{A_{j}^{l}}\mathrm{d}X_{s}^{t,x;\tilde{v}^{l,j},\tilde{\xi}^{l,j}}
& = & \mathbf{1}_{A_{j}^{l}}b\left( s,X_{s}^{t,x;\tilde{v}^{l,j},\tilde{\xi}%
^{l,j}},\tilde{v}_{s}^{l,j}\right) \mathrm{d}s \\
&  & +\mathbf{1}_{A_{j}^{l}}\sigma \left( s,X_{s}^{t,x;\tilde{v}^{l,j},%
\tilde{\xi}^{l,j}},\tilde{v}_{s}^{l,j}\right) \mathrm{d}W_{s}+\mathbf{1}%
_{A_{j}^{l}}G\mathrm{d}\tilde{\xi}_{s}^{l,j}, \\
\mathbf{1}_{A_{j}^{l}}\mathrm{d}Y_{s}^{t,x;\tilde{v}^{l,j},\tilde{\xi}^{l,j}}
& = & -\mathbf{1}_{A_{j}^{l}}f\left( t,X_{s}^{t,x;\tilde{v}^{l,j},\tilde{\xi}%
^{l,j}},Y_{s}^{t,x;\tilde{v}^{l,j},\tilde{\xi}^{l,j}},Z_{s}^{t,x;\tilde{v}%
^{l,j},\tilde{\xi}^{l,j}},\tilde{v}_{s}^{l,j}\right) \mathrm{d}s \\
&  & +\mathbf{1}_{A_{j}^{l}}Z_{s}^{t,x;\tilde{v}^{l,j},\tilde{\xi}^{l,j}}%
\mathrm{d}W_{s}-K\mathbf{1}_{A_{j}^{l}}\mathrm{d}\tilde{\xi}_{s}^{l,j}, \\
X_{t}^{t,x;\tilde{v}^{l,j},\tilde{\xi}^{l,j}} & = & x\mathbf{1}_{A_{j}^{l}},%
\text{ }Y_{T}^{t,x;\tilde{v}^{l,j},\tilde{\xi}^{l,j}}=\mathbf{1}%
_{A_{j}^{l}}\Phi \left( X_{T}^{t,x;\tilde{v}^{l,j},\tilde{\xi}^{l,j}}\right)
,0\leq t\leq s\leq T.%
\end{array}%
\right.
\end{equation*}%
It is easy to derive that
\begin{equation*}
\left\{
\begin{array}{lll}
\mathrm{d}\sum_{j=0}^{N_{l}-1}\mathbf{1}_{A_{j}^{l}}X_{s}^{t,x;\tilde{v}%
^{l,j},\tilde{\xi}^{l,j}} & = & b\left( s,\sum_{j=0}^{N_{l}-1}\mathbf{1}%
_{A_{j}^{l}}X_{s}^{t,x;\tilde{v}^{l,j},\tilde{\xi}^{l,j}},%
\sum_{j=0}^{N_{l}-1}\mathbf{1}_{A_{j}^{l}}\tilde{v}_{s}^{l,j}\right) \mathrm{%
d}s \\
&  & +\sigma \left( s,\sum_{j=0}^{N_{l}-1}\mathbf{1}_{A_{j}^{l}}X_{s}^{t,x;%
\tilde{v}^{l,j},\tilde{\xi}^{l,j}},\sum_{j=0}^{N_{l}-1}\mathbf{1}_{A_{j}^{l}}%
\tilde{v}_{s}^{l,j}\right) \mathrm{d}W_{s} \\
&  & +G\sum_{j=0}^{N_{l}-1}\mathbf{1}_{A_{j}^{l}}\mathrm{d}\tilde{\xi}%
_{s}^{l,j}, \\
\mathrm{d}\sum_{j=0}^{N_{l}-1}\mathbf{1}_{A_{j}^{l}}Y_{s}^{t,x;\tilde{v}%
^{l,j},\tilde{\xi}^{l,j}} & = & -f\Big (t,\sum_{j=0}^{N_{l}-1}\mathbf{1}%
_{A_{j}^{l}}X_{s}^{t,x;\tilde{v}^{l,j},\tilde{\xi}^{l,j}},%
\sum_{j=0}^{N_{l}-1}\mathbf{1}_{A_{j}^{l}}Y_{s}^{t,x;\tilde{v}^{l,j},\tilde{%
\xi}^{l,j}}, \\
&  & \sum_{j=0}^{N_{l}-1}\mathbf{1}_{A_{j}^{l}}Z_{s}^{t,x;\tilde{v}^{l,j},%
\tilde{\xi}^{l,j}},\sum_{j=0}^{N_{l}-1}\mathbf{1}_{A_{j}^{l}}\tilde{v}%
_{s}^{l,j}\Big )\mathrm{d}s \\
&  & +\sum_{j=0}^{N_{l}-1}\mathbf{1}_{A_{j}^{l}}Z_{s}^{t,x;\tilde{v}^{l,j},%
\tilde{\xi}^{l,j}}\mathrm{d}W_{s}-K\sum_{j=0}^{N_{l}-1}\mathbf{1}_{A_{j}^{l}}%
\mathrm{d}\tilde{\xi}_{s}^{l,j}, \\
\sum_{j=0}^{N_{l}-1}\mathbf{1}_{A_{j}^{l}}X_{t}^{t,x;\tilde{v}^{l,j},\tilde{%
\xi}^{l,j}} & = & \sum_{j=0}^{N_{l}-1}\mathbf{1}_{A_{j}^{l}}x, \\
\sum_{j=0}^{N_{l}-1}\mathbf{1}_{A_{j}^{l}}Y_{T}^{t,x;\tilde{v}^{l,j},\tilde{%
\xi}^{l,j}} & = & \Phi \left( \sum_{j=0}^{N_{l}-1}\mathbf{1}%
_{A_{j}^{l}}X_{T}^{t,x;\tilde{v}^{l,j},\tilde{\xi}^{l,j}}\right) ,0\leq
t\leq s\leq T.%
\end{array}%
\right.
\end{equation*}%
Combining the uniqueness of FBSDEs, we claim that%
\begin{eqnarray*}
X_{s}^{t,x;\tilde{v}^{l},\tilde{\xi}^{l}} &=&\sum_{j=0}^{N_{l}-1}\mathbf{1}%
_{A_{j}^{l}}X_{s}^{t,x;\tilde{v}^{l,j},\tilde{\xi}^{l,j}}, \\
Y_{s}^{t,x;\tilde{v}^{l},\tilde{\xi}^{l}} &=&\sum_{j=0}^{N_{l}-1}\mathbf{1}%
_{A_{j}^{l}}Y_{s}^{t,x;\tilde{v}^{l,j},\tilde{\xi}^{l,j}}, \\
Z_{s}^{t,x;\tilde{v}^{l},\tilde{\xi}^{l}} &=&\sum_{j=0}^{N_{l}-1}\mathbf{1}%
_{A_{j}^{l}}Z_{s}^{t,x;\tilde{v}^{l,j},\tilde{\xi}^{l,j}},\text{ }s\in \left[
t,T\right] .
\end{eqnarray*}%
From $J(t,x;\tilde{v}^{l,j},\tilde{\xi}^{l,j}):=\left. Y_{s}^{t,x;\tilde{v}%
^{l,j},\tilde{\xi}^{l,j}}\right\vert _{s=t},$ we conclude that%
\begin{eqnarray*}
J(t,x;\tilde{v}^{l},\tilde{\xi}^{l}) &=&\sum_{j=0}^{N_{l}-1}\mathbf{1}%
_{A_{j}^{l}}\left. Y_{s}^{t,x;\tilde{v}^{l,j},\tilde{\xi}^{l,j}}\right\vert
_{s=t} \\
&=&\sum_{j=0}^{N_{l}-1}\mathbf{1}_{A_{j}^{l}}J(t,x;\tilde{v}^{l,j},\tilde{\xi%
}^{l,j}).
\end{eqnarray*}%
Since $Y_{s}^{t,x;\tilde{v}^{l,j},\tilde{\xi}^{l,j}}$ is $\mathscr{F}%
_{s}^{t} $-adapted, we know that $J(t,x;\tilde{v}^{l,j},\tilde{\xi}%
^{l,j})\in \mathbb{R}$. Moreover,
\begin{equation*}
u\left( t,x\right) \leq \min_{0\leq j\leq N_{l-1}}J\left( t,x;\tilde{v}%
^{l,j},\tilde{\xi}^{l,j}\right) \leq J(t,x;\tilde{v}^{l},\tilde{\xi}^{l}).
\end{equation*}%
Assume that $j=j_{l}$ such that
\begin{equation*}
\min_{0\leq j\leq N_{l-1}}J\left( t,x;\tilde{v}^{l,j},\tilde{\xi}%
^{l,j}\right) =J\left( t,x;\tilde{v}^{l,j_{l}},\tilde{\xi}^{l,j_{l}}\right) .
\end{equation*}%
By (\ref{coner1}), we can get
\begin{equation*}
u\left( t,x\right) =\lim_{l\rightarrow +\infty }J\left( t,x;\tilde{v}%
^{l,j_{l}},\tilde{\xi}^{l,j_{l}}\right) \in \mathbb{R}\text{.}
\end{equation*}%
The proof is completed. ~\hfill $\Box $

\begin{lemma}
\label{adcontrol}For any $\varepsilon >0,$ there exists a control $\left(
v,\xi \right) \in \mathcal{U}$, such that
\begin{equation*}
u\left( t,x\right) \geq Y_{s}^{t,x;v,\xi }-\varepsilon ,
\end{equation*}%
where $u\left( t,x\right) $ is defined in \emph{(\ref{value})}.
\end{lemma}

\paragraph{Proof}

From (\ref{coner1}), we can find out an $\mathscr{F}_{s}^{t}$ adapted
control pair $\left( v^{i},\xi ^{i}\right) $ such that $u\left( t,x\right)
\geq Y_{s}^{t,x;v^{i},\xi ^{i}}-\varepsilon ,$ a.s.. We denote that $\left(
v,\xi \right) =\sum_{i=1}^{N}\left( v^{i},\xi ^{i}\right) \mathbf{1}%
_{A_{i}}, $ where $A_{i}$ is a partition of $\left( \Omega ,\mathscr{F}%
_{t}\right) .$ Next, we proceed that%
\begin{eqnarray*}
Y_{s}^{t,x;v,\xi } &=&\sum_{i=1}^{N}Y_{s}^{t,x;v^{i},\xi ^{i}}\mathbf{1}%
_{A_{i}} \\
&\leq &\sum_{i=1}^{N}\left( u\left( t,x\right) +\varepsilon \right) \mathbf{1%
}_{A_{i}} \\
&=&u\left( t,x\right) +\varepsilon ,
\end{eqnarray*}%
from which we get the desired result. ~\hfill $\Box $

\begin{lemma}
\label{DPP}Under assumptions \emph{(A1)--(A3)}, the value function $%
u\left( t,x\right) $ defined in \emph{(\ref{value})} obeys the DPP: For any $%
0\leq t\leq t+\delta \leq T,$ $x\in \mathbb{R}^{n},$%
\begin{equation*}
u\left( t,x\right) =ess\inf_{\left( v,\xi \right) \in \mathcal{U}}\mathcal{G}%
_{t,t+\delta }^{t,x,v,\xi }\left[ u\left( t+\delta ,X_{t+\delta }^{t,x,v,\xi
}\right) \right] .
\end{equation*}
\end{lemma}

\paragraph{Proof}

First we have
\begin{eqnarray*}
u\left( t,x\right) &=&ess\inf_{\left( v,\xi \right) \in \mathcal{U}}\mathcal{%
G}_{t,T}^{t,x,v,\xi }\left[ \Phi \left( X_{T}^{t,x,v,\xi }\right) \right] \\
&=&ess\inf_{\left( v,\xi \right) \in \mathcal{U}}\mathcal{G}_{t,t+\delta
}^{t,x,v,\xi }\left[ Y_{t+\delta }^{t+\delta ,X_{t+\delta }^{t,x,v},v,\xi }%
\right] .
\end{eqnarray*}%
On the one hand, by comparison theorem (lemma \ref{comp}), we have%
\begin{equation*}
u\left( t,x\right) \geq ess\inf_{\left( v,\xi \right) \in \mathcal{U}}%
\mathcal{G}_{t,T}^{t,x,v,\xi }\left[ u\left( t+\delta ,X_{t+\delta
}^{t,x,v,\xi }\right) \right] .
\end{equation*}%
On the other hand, for any $\varepsilon >0,$ there exists an admissible
control pair $\left( \bar{v},\bar{\xi}\right) \in \mathcal{U}$ such that%
\begin{equation*}
u\left( t+\delta ,X_{t+\delta }^{t,x;\bar{v},\bar{\xi}}\right) \geq
Y_{t+\delta }^{t+\delta ,X_{t+\delta }^{t,x;\bar{v},\bar{\xi}};\bar{v},\bar{%
\xi}}-\varepsilon ,\text{ a.s.}
\end{equation*}%
Then, for any $\left( v,\xi \right) \in \mathcal{U}$, we define
\begin{equation*}
\left( \tilde{v},\tilde{\xi}\right) _{s}=\left\{
\begin{array}{c}
\left( v,\xi \right) _{s},\text{ }s\in \left[ 0,t+\delta \right] , \\
\left( \bar{v},\bar{\xi}\right) _{s},\text{ }s\in \left[ t+\delta ,T\right] .%
\end{array}%
\right.
\end{equation*}%
Clearly, $\left( \tilde{v},\tilde{\xi}\right) \in \mathcal{U}$. By the
comparison theorem (lemma \ref{comp}) again, we obtain%
\begin{eqnarray*}
u\left( t,x\right) &\leq &\mathcal{G}_{t,t+\delta }^{t,x;\tilde{v},\tilde{\xi%
}}\left[ Y_{t+\delta }^{t+\delta ,X_{t+\delta }^{t,x;\tilde{v},\tilde{\xi}};%
\tilde{v},\tilde{\xi}}\right] \\
&=&\mathcal{G}_{t,t+\delta }^{t,x;v,\xi }\left[ Y_{t+\delta }^{t+\delta
,X_{t+\delta }^{t,x;v,\xi };\bar{v},\bar{\xi}}\right] \\
&\leq &\mathcal{G}_{t,t+\delta }^{t,x;v,\xi }\left[ u\left( t+\delta
,X_{t+\delta }^{t,x;v,\xi }\right) +\varepsilon \right] \\
&\leq &\mathcal{G}_{t,t+\delta }^{t,x;v,\xi }\left[ u\left( t+\delta
,X_{t+\delta }^{t,x;v,\xi }\right) \right] +C\varepsilon ,
\end{eqnarray*}%
where $C$ is independent of control pair. We get the DPP for $u\left(
t,x\right) $ since the arbitrary of $\varepsilon .$ ~\hfill $\Box $

\begin{lemma}
\label{LH}Under assumptions \emph{(A1)--(A3)}, the value function $%
u\left( t,x\right) $ defined in \emph{(\ref{value})} admits%
\begin{equation*}
\left\vert u\left( t,x\right) -u\left( t^{\prime },x^{\prime }\right)
\right\vert \leq C\left( \left\vert x-x^{\prime }\right\vert +\left\vert
t-t^{\prime }\right\vert ^{\frac{1}{2}}\right) ,\text{ }0\leq t,t^{\prime
}\leq T,\text{ }x,x^{\prime }\in \mathbb{R}^{n}.
\end{equation*}
\end{lemma}

\paragraph{Proof}

For any $\left( t,x\right) \in \left[ 0,T\right] \times \mathbb{R}^{n}$ and $%
\forall \delta >0,$ from Lemma \ref{DPP} and Lemma \ref{adcontrol}, there
exists, for any $\varepsilon >0,$ an admissible control pair $\left(
v^{\varepsilon },\xi ^{\varepsilon }\right) $ such that
\begin{equation}
\mathcal{G}_{t,t+\delta }^{t,x;v^{\varepsilon },\xi ^{\varepsilon }}\left[
u\left( t+\delta ,X_{t+\delta }^{t,x;v^{\varepsilon },\xi ^{\varepsilon
}}\right) \right] -\varepsilon \leq u\left( t,x\right) \leq \mathcal{G}%
_{t,t+\delta }^{t,x;v^{\varepsilon },\xi ^{\varepsilon }}\left[ u\left(
t+\delta ,X_{t+\delta }^{t,x;v^{\varepsilon },\xi ^{\varepsilon }}\right) %
\right] .  \label{holder}
\end{equation}%
From (\ref{holder}), we focus on
\begin{equation*}
u\left( t,x\right) -u\left( t+\delta ,x\right) \geq \Theta _{\delta
}^{1}+\Theta _{\delta }^{1}+\varepsilon ,
\end{equation*}%
where
\begin{eqnarray*}
\Theta _{\delta }^{1} &=&\mathcal{G}_{t,t+\delta }^{t,x;v^{\varepsilon },\xi
^{\varepsilon }}\left[ u\left( t+\delta ,X_{t+\delta }^{t,x;v^{\varepsilon
},\xi ^{\varepsilon }}\right) \right] -\mathcal{G}_{t,t+\delta
}^{t,x;v^{\varepsilon },\xi ^{\varepsilon }}\left[ u\left( t+\delta
,x\right) \right] , \\
\Theta _{\delta }^{2} &=&\mathcal{G}_{t,t+\delta }^{t,x;v^{\varepsilon },\xi
^{\varepsilon }}\left[ u\left( t+\delta ,x\right) \right] -u\left( t+\delta
,x\right) .
\end{eqnarray*}%
From Lemma \ref{LH}, we get
\begin{eqnarray*}
\left\vert \Theta _{\delta }^{1}\right\vert &\leq &C\mathbb{E}^{\mathscr{F}%
_{t}}\left[ \left\vert u\left( t+\delta ,X_{t+\delta }^{t,x;v^{\varepsilon
},\xi ^{\varepsilon }}\right) -u\left( t+\delta ,x\right) \right\vert ^{2}%
\right] ^{\frac{1}{2}} \\
&\leq &C\mathbb{E}^{\mathscr{F}_{t}}\left[ \left\vert X_{t+\delta
}^{t,x;v^{\varepsilon },\xi ^{\varepsilon }}-x\right\vert ^{2}\right] ^{%
\frac{1}{2}} \\
&\leq &C\left( 1+\left\vert x\right\vert ^{2}\right) \delta .
\end{eqnarray*}%
Hence, we deduce that
\begin{equation}
\Theta _{\delta }^{1}\geq -C\left( 1+\left\vert x\right\vert \right) \delta
^{\frac{1}{2}},  \label{esttime1}
\end{equation}%
since $\sqrt{\left( 1+\left\vert x\right\vert ^{2}\right) }\leq 1+\left\vert
x\right\vert ,$ $\forall x\in \mathbb{R}^{n}.$ Next we deal with $\Theta
_{\delta }^{2}.$ From the definition of DDP, we know that $\Theta _{\delta
}^{2}$ can be rewritten as
\begin{eqnarray*}
\Theta _{\delta }^{2} &=&\mathbb{E}^{\mathscr{F}_{t}}\left[ u\left( t+\delta
,x\right) +\int_{t}^{t+\delta }f\left( s,X_{s}^{t,x;v^{\varepsilon },\xi
^{\varepsilon }},Y_{s}^{t,x;v^{\varepsilon },\xi ^{\varepsilon
}},Z_{s}^{t,x;v^{\varepsilon },\xi ^{\varepsilon }},v_{t}^{\varepsilon
}\right) \mathrm{d}s\right. \\
&&\left. -\int_{t}^{t+\delta }Z_{s}^{t,x;v^{\varepsilon },\xi ^{\varepsilon
}}\mathrm{d}W_{s}+\int_{t}^{t+\delta }K\mathrm{d}\xi _{s}^{\varepsilon }%
\right] -u\left( t+\delta ,x\right) \\
&=&\mathbb{E}^{\mathscr{F}_{t}}\left[ \int_{t}^{t+\delta }f\left(
s,X_{s}^{t,x;v^{\varepsilon },\xi ^{\varepsilon }},Y_{s}^{t,x;v^{\varepsilon
},\xi ^{\varepsilon }},Z_{s}^{t,x;v^{\varepsilon },\xi ^{\varepsilon
}},v_{t}^{\varepsilon }\right) \mathrm{d}s+\int_{t}^{t+\delta }K\mathrm{d}%
\xi _{s}^{\varepsilon }\right] .
\end{eqnarray*}%
Employing Schwartz inequality, we have, for certain constant $C$ not
depending on $t$ and $\delta ,$%
\begin{eqnarray*}
\left\vert \Theta _{\delta }^{2}\right\vert &\leq &C\delta ^{\frac{1}{2}}%
\mathbb{E}^{\mathscr{F}_{t}}\Bigg [\int_{t}^{t+\delta }\Big |f\left(
s,X_{s}^{t,x;v^{\varepsilon },\xi ^{\varepsilon }},0,0,v_{t}^{\varepsilon
}\right) +\left\vert Y_{s}^{t,x;v^{\varepsilon },\xi ^{\varepsilon
}}\right\vert \\
&&+\left\vert Z_{s}^{t,x;v^{\varepsilon },\xi ^{\varepsilon }}\right\vert %
\Big |\mathrm{d}s+\left\vert \xi _{T}^{\varepsilon }\right\vert ^{2}\Bigg ]
\\
&\leq &C\delta ^{\frac{1}{2}}\mathbb{E}^{\mathscr{F}_{t}}\Bigg [%
\int_{t}^{t+\delta }\left( 1+\left\vert X_{s}^{t,x;v^{\varepsilon },\xi
^{\varepsilon }}\right\vert ^{2}+\left\vert Y_{s}^{t,x;v^{\varepsilon },\xi
^{\varepsilon }}\right\vert ^{2}+\left\vert Z_{s}^{t,x;v^{\varepsilon },\xi
^{\varepsilon }}\right\vert ^{2}\right) \mathrm{d}s \\
&&+\left\vert \xi _{T}^{\varepsilon }\right\vert ^{2}\Bigg ] \\
&\leq &C\left( 1+\left\vert x\right\vert ^{2}\right) \delta ^{\frac{1}{2}},
\end{eqnarray*}%
that is
\begin{equation}
\Theta _{\delta }^{2}\geq -C\left( 1+\left\vert x\right\vert ^{2}\right)
\delta ^{\frac{1}{2}}.  \label{estime2}
\end{equation}%
Therefore, combing (\ref{esttime1}) and (\ref{estime2}), we have%
\begin{equation*}
u\left( t,x\right) -u\left( t+\delta ,x\right) \geq -C\left( 1+\left\vert
x\right\vert ^{2}\right) \delta ^{\frac{1}{2}}-\varepsilon .
\end{equation*}%
Note that the arbitrary of $\varepsilon >0,$ we immitigably get
\begin{equation*}
u\left( t,x\right) -u\left( t+\delta ,x\right) \geq -C\left( 1+\left\vert
x\right\vert ^{2}\right) \delta ^{\frac{1}{2}}.
\end{equation*}%
Repeating the above approach, we can get the desired result basing on the
second inequality of (\ref{holder}). ~\hfill $\Box $

We now study the original optimal control problem via DPP (Lemma \ref{DPP}).
Postulate that $u\left( t,x\right) \in C^{1,2}\left( \left[ 0,T\right]
\times \mathbb{R}^{n};\mathbb{R}\right) .$ Indeed, this hypothesis is really
strong, but the main destination here is to derive the so called H-J-B
inequality.

Define
\begin{eqnarray*}
\mathcal{L}\left( t,x,v\right) \Psi  &=&\frac{1}{2}\text{\textrm{Tr}}\left(
\sigma \sigma ^{\ast }\left( t,x,v\right) D^{2}\Psi \right) +\left\langle
D\Psi ,b\left( t,x,v\right) \right\rangle , \\
\left( t,x,v\right)  &\in &\left[ 0,T\right] \times \mathbb{R}^{n}\times U,%
\text{ }\Psi \in C^{1,2}\left( \left[ 0,T\right] \times \mathbb{R}%
^{n}\right) .
\end{eqnarray*}%
For any $\Psi \in C^{1,2}\left( \left[ 0,T\right] \times \mathbb{R}^{n};%
\mathbb{R}\right) $, by virtue of Dol\'{e}ans--Dade--Meyer formula (see \cite%
{HS2,CH, W}), we have%
\begin{eqnarray}
\Psi \left( s,X_{s}\right)  &=&\Psi \left( t,x\right) +\int_{t}^{s}\Psi
_{t}\left( r,X_{r}\right) +\mathcal{L}\left( r,X_{r},v\right) \Psi \left(
r,X_{r}\right) \mathrm{d}r  \notag \\
&&+\int_{t}^{s}\Psi _{x}\left( r,X_{r}\right) \sigma \left(
r,X_{r},v_{r}\right) \mathrm{d}W_{r}+\int_{t}^{s}\Psi _{x}\left(
r,X_{r}\right) G\mathrm{d}\xi _{r}  \notag \\
&&+\sum_{t\leq r\leq s}\left\{ \Psi \left( r,X_{r+}\right) -\Psi \left(
r,X_{r}\right) -\Psi _{x}\left( r,X_{r}\right) \Delta X_{r}\right\} .
\label{ito1}
\end{eqnarray}%
Taking expectation, we get
\begin{eqnarray}
&&\mathbb{E}\left[ \Psi \left( s,X_{s}\right) \right]   \notag \\
&=&\Psi \left( t,X_{t}\right) +\mathbb{E}\left[ \int_{t}^{s}\Psi _{x}\left(
r,X_{r}\right) G\mathrm{d}\xi _{r}\right]   \notag \\
&&+\mathbb{E}\left[ \int_{t}^{s}\Psi _{t}\left( r,X_{r}\right) +\mathcal{L}%
\left( t,x,v\right) \Psi \left( r,X_{r}\right) \mathrm{d}r\right]   \notag \\
&&+\mathbb{E}\left[ \sum_{t\leq r\leq s}\left\{ \Psi \left( r,X_{r+}\right)
-\Psi \left( r,X_{r}\right) -\Psi _{x}\left( r,X_{r}\right) \Delta
X_{r}\right\} \right] .  \label{ito2}
\end{eqnarray}%
Now according to the DPP (Lemma \ref{DPP}), it follows%
\begin{eqnarray}
u\left( t,x\right)  &=&\inf_{v\left( \cdot \right) \times \xi \left( \cdot
\right) \in \mathcal{U}}\mathbb{E}\Bigg [\int_{t}^{s}f\left(
r,X_{r}^{t,x;v,\xi },Y_{r}^{t,x;v,\xi },Z_{r}^{t,x;v,\xi },v_{r}\right)
\mathrm{d}r  \notag \\
&&-\int_{t}^{s}Z_{r}^{t,x;v,\xi }\mathrm{d}W_{r}+\int_{t}^{s}K\mathrm{d}\xi
_{r}+u\left( s,X_{s}^{t,x;v,\xi }\right) \Bigg ].  \label{d2}
\end{eqnarray}%
Hence, by It\^{o}'s formula (\ref{ito2}) again, (\ref{d2}) yields%
\begin{eqnarray*}
0 &=&\inf_{v\left( \cdot \right) \times \xi \left( \cdot \right) \in
\mathcal{U}}\mathbb{E}\Bigg [\int_{t}^{s}\Big [u_{t}\left(
r,X_{r}^{t,x;v,\xi }\right) +\mathcal{L}\left( r,X_{r}^{t,x;v,\xi
},v_{r}\right) u\left( r,X_{r}^{t,x;v,\xi }\right)  \\
&&+f\left( r,X_{r}^{t,x;v,\xi },Y_{r}^{t,x;v,\xi },Z_{r}^{t,x;v,\xi
},v_{r}\right) \Big ]\mathrm{d}r \\
&&+\int_{t}^{s}\left( u_{x}^{\top }\left( r,X_{r}^{t,x;v,\xi }\right)
G+K\right) \mathrm{d}\xi _{r} \\
&&+\int_{t}^{s}\left( u_{x}^{\top }\left( r,X_{r}\right) \sigma \left(
r,X_{r},v_{r}\right) -Z_{r}^{t,x;v,\xi }\right) \mathrm{d}W_{r}
\end{eqnarray*}

\begin{eqnarray}
&&+\sum_{t\leq r\leq s}\Big [u\left( r,X_{r+}^{t,x;v,\xi }\right) -u\left(
r,X_{r}^{t,x;v,\xi }\right)  \notag \\
&&-u_{x}\left( r,X_{r}^{t,x;v,\xi }\right) \Delta X_{r}^{t,x;v,\xi }\Big ]%
\Bigg ].  \label{dy1}
\end{eqnarray}%
From (\ref{dy1}), it is possible to give the following definition of
classical solution.

\begin{definition}[Classical Solution]
\label{df1}Consider a function $V\left( t,x\right) \in C^{1,2}\left( \left[
0,T\right] \times \mathbb{R}^{n};\mathbb{R}\right) $ and define that
\begin{equation*}
\mathcal{D}_{t}\left( V\right) :=\left\{ x\in \mathbb{R}^{n}:V_{x}^{\top
}\left( t,x\right) G_{t}+K>0\right\} .
\end{equation*}%
We can say $V$ is a classical solution of the dynamic programming inequality
\emph{(\ref{dy1})}, if for every $\left( t,x\right) \in \left[ 0,T\right]
\times \mathcal{D}_{t}\left( V\right) ,$ $V$ admits%
\begin{equation*}
\frac{\partial }{\partial t}V\left( t,x\right) +\min_{v\in U}\left\{
\mathcal{L}\left( t,x,v\right) V\left( t,x\right) +f\left( t,x,V\left(
t,x\right) ,\nabla V\left( t,x\right) \sigma \left( t,x,v\right) ,v\right)
\right\} =0,
\end{equation*}%
and
\begin{equation*}
V\left( T,x\right) =\Phi \left( x\right) ,\text{ }x\in \mathbb{R}^{n}.
\end{equation*}%
Moreover, for any $\left( t,x,v\right) \in \left[ 0,T\right] \times \mathbb{R%
}^{n}\times U,$
\begin{equation}
V_{x}^{\top }\left( t,x\right) G+K\geq 0,  \label{b1}
\end{equation}%
and%
\begin{equation*}
\frac{\partial }{\partial t}V\left( t,x\right) +\mathcal{L}\left(
t,x,v\right) V\left( t,x\right) +f\left( t,x,V\left( t,x\right) ,\nabla
V\left( t,x\right) \sigma \left( t,x,v\right) ,v\right) \geq 0.
\end{equation*}
\end{definition}

\begin{remark}
In the Definition (\ref{df1}), $V_{x}^{\top }\left( t,x\right) G+K\geq 0$
means $\left( V_{x}^{\top }\left( t,x\right) G\right) ^{i}+K^{i}\geq 0$, $%
i=1,2,\ldots ,m,$ where $\left( \cdot \right) ^{i}$ denotes the $i$th
coordinate of the point in $\mathbb{R}^{m}.$ Now we introduce the following:
\begin{eqnarray}
&&0=\min \Big \{\frac{\partial }{\partial t}V\left( t,x\right)  \notag \\
&&+\min_{v\in U}\left\{ \mathcal{L}\left( t,x,v\right) V\left( t,x\right)
+f\left( t,x,V\left( t,x\right) ,\nabla V\left( t,x\right) \sigma \left(
t,x,v\right) ,v\right) \right\} ,  \notag \\
&&V_{x}^{\top }\left( t,x\right) G+K\Big \}.  \label{dy2}
\end{eqnarray}
\end{remark}

%

We shall prove that, under certain assumptions, the classical solution $V$ to
H-J-B inequality (\ref{dy2}) is just the optimal value function with optimal
control pair.

\begin{theorem}[Verification Theorem]
\label{t2} Suppose that $V$ is a classical solution of the H-J-B inequality
\emph{(\ref{dy2})} such that for some $l>1,$ $\left\vert V\left( t,x\right)
\right\vert \leq C\left( 1+\left\vert x\right\vert ^{l}\right) .$ Then for
any $\left[ 0,T\right] \times \mathbb{R}^{n}$, $\left( v,\xi \right) \in
\mathcal{U}:$%
\begin{equation*}
V\left( t,x\right) \leq J\left( t,x,v,\xi \right) .
\end{equation*}%
Furthermore, if there exists $\left( \hat{v},\hat{\xi}\right) \in \mathcal{U}
$ such that
\begin{equation}
P\left\{ \left( r,X_{r}^{t,x;\hat{v},\hat{\xi}}\right) \in \mathcal{D}%
_{r}\left( V\right) ,\text{ }0\leq r\leq T\right\} =1,  \label{v11}
\end{equation}%
\begin{equation}
P\left\{ \int_{\left[ t,T\right] }\left[ V_{x}^{\top }\left( r,x\right) G+K%
\right] \mathrm{d}r=0\right\} =1,  \label{v22}
\end{equation}%
\begin{eqnarray}
1 &=&P\Bigg \{\left( s,X_{s+}^{t,x;\hat{v},\hat{\xi}}\right) \in \mathcal{D}%
_{r}\left( V\right) ,\text{ }t\leq s\leq T:  \notag \\
\hat{v}_{s} &\in &\min_{v\in U}\Big [V_{t}\left( s,X_{s}^{t,x;v,\xi }\right)
+\mathcal{L}\left( s,X_{s+}^{t,x;\hat{v},\hat{\xi}},v\right) V\left(
s,X_{s+}^{t,x;\hat{v},\hat{\xi}}\right)  \notag \\
&&+f\Big (s,X_{s+}^{t,x;\hat{v},\hat{\xi}},V\left( t,X_{s+}^{t,x;\hat{v},%
\hat{\xi}}\right) ,  \notag \\
&&\nabla V\left( s,X_{s+}^{t,x;\hat{v},\hat{\xi}}\right) \sigma \left(
s,X_{s+}^{t,x;\hat{v},\hat{\xi}},v\right) ,v\Big )\Big ]\Bigg \}  \label{v33}
\end{eqnarray}%
and%
\begin{equation}
P\left\{ V\left( s,X_{s}^{t,x;\hat{v},\hat{\xi}}\right) =V\left(
s,X_{s+}^{t,x;\hat{v},\hat{\xi}}\right) +K\Delta \hat{\xi}_{s},\text{ }t\leq
s\leq T\right\} =1.  \label{v44}
\end{equation}%
Then%
\begin{equation}
V\left( t,x\right) =J\left( t,x;\hat{v}\left( \cdot \right) ,\hat{\xi}\left(
\cdot \right) \right) .  \label{v55}
\end{equation}
\end{theorem}

\paragraph{Proof}

For any $\left( t,x\right) \in \left[ 0,T\right] \times \mathbb{R}^{n}$ and $%
\left( v,\xi \right) \in \mathcal{U},$ applying Dol\'{e}ans--Dade--Meyer
formula to $V,$ we have%
\begin{eqnarray}
V\left( t,x\right) &=&\mathbb{E}\left[ \Phi \left( X_{T}^{t,x;v,\xi }\right) %
\right] -\mathbb{E}\Bigg [\int_{t}^{s}\left( r,X_{r}^{t,x;v,\xi }\right)
\notag \\
&&+\mathcal{L}\left( r,X_{r}^{t,x;v,\xi },v\right) V\left(
r,X_{r}^{t,x;v,\xi }\right) \mathrm{d}r\Bigg ]  \notag \\
&&-\mathbb{E}\left[ \int_{t}^{s}V_{x}\left( r,X_{r}^{t,x;v,\xi }\right) G%
\mathrm{d}\xi _{r}\right]  \notag \\
&&-\mathbb{E}\Bigg [\sum_{t\leq r\leq s}\left[ V\left( r,X_{r+}^{t,x;v,\xi
}\right) -V\left( r,X_{r}^{t,x;v,\xi }\right) \right.  \notag \\
&&\left. -V_{x}\left( r,X_{r}^{t,x;v,\xi }\right) \Delta X_{r}^{t,x;v,\xi }
\right] \Bigg ].  \label{v1}
\end{eqnarray}%
Note that $\Delta X_{r}^{t,x;v,\xi }=G\Delta \xi _{r}$, and $%
X_{r+}^{t,x;v,\xi }=X_{r}^{t,x;v,\xi }+\Delta X_{r}^{t,x;v,\xi
}=X_{r}^{t,x;v,\xi }+G\Delta \xi _{r}.$ Thus%
\begin{eqnarray}
&&-\mathbb{E}\left[ \int_{t}^{s}V_{x}\left( r,X_{r}^{t,x;v,\xi }\right) G%
\mathrm{d}\xi _{r}\right] +\mathbb{E}\left[ V_{x}\left( r,X_{r}^{t,x;v,\xi
}\right) \Delta X_{r}^{t,x;v,\xi }\right]  \notag \\
&=&-\mathbb{E}\left[ \int_{t}^{s}V_{x}\left( r,X_{r}^{t,x;v,\xi }\right) G%
\mathrm{d}\xi _{r}^{c}\right]  \notag \\
&\leq &\mathbb{E}\left[ \int_{t}^{s}K\mathrm{d}\xi _{r}^{c}\right] .
\label{v2}
\end{eqnarray}%
We now deal the term
\begin{eqnarray}
&&-\mathbb{E}\left[ \sum_{t\leq r\leq s}\left\{ V\left( r,X_{r+}^{t,x;v,\xi
}\right) -V\left( r,X_{r}^{t,x;v,\xi }\right) \right\} \right]  \notag \\
&=&-\mathbb{E}\left[ \sum_{t\leq r\leq s}\left\{ \int_{0}^{1}V_{x}\left(
r,X_{r}^{t,x;v,\xi }+\theta \Delta X_{r}^{t,x;v,\xi }\right) G\Delta \xi _{r}%
\mathrm{d}\theta \right\} \right]  \notag \\
&\leq &\mathbb{E}\left[ K\Delta \xi _{r}\right] .  \label{v3}
\end{eqnarray}%
Combining (\ref{v2}) and (\ref{v3}), we have%
\begin{eqnarray*}
V\left( t,x\right) &\leq &\mathbb{E}\Bigg [\Phi \left( X_{T}^{t,x;v,\xi
}\right) +\int_{t}^{s}f\left( r,X_{r}^{t,x;v,\xi },Y_{r}^{t,x;v,\xi
},Z_{r}^{t,x;v,\xi },v_{r}\right) \mathrm{d}r \\
&&+\int_{t}^{s}K\mathrm{d}\xi _{r}\Bigg ].
\end{eqnarray*}%
Now we take $\left( \hat{v},\hat{\xi}\right) \in \mathcal{U}$ such that (\ref%
{v11})-(\ref{v44}) hold. Then repeating above method, we have%
\begin{eqnarray*}
V\left( t,x\right) &\leq &\mathbb{E}\Bigg [\Phi \left( X_{T}^{t,x;\hat{v},%
\hat{\xi}}\right) +\int_{t}^{s}f\left( r,X_{r}^{t,x;\hat{v},\hat{\xi}%
},Y_{r}^{t,x;\hat{v},\hat{\xi}},Z_{r}^{t,x;\hat{v},\hat{\xi}},\hat{v}%
_{r}\right) \mathrm{d}r \\
&&+\int_{t}^{s}K\mathrm{d}\hat{\xi}_{r}\Bigg ].
\end{eqnarray*}%
We thus complete the proof. \hfill $\Box $

\begin{remark}
As observed in Theorem \ref{t2}, it provides a possible way to check a pair
of potential singular controls is optimal or not. Anyhow, note that Section
6 in \cite{W} is devoted to this issue via constructing an optimal control.
In this sense, our verification theorem seems to be rather weak, since we do
not know whether $(\hat{v},\hat{\xi})$ exists or not. Regarding to this
topic, we shall discus in Remark \ref{existence} below. Verification theorem
under non-smooth situation, namely, viscosity approach can be found in Zhang
\cite{zlq}.
\end{remark}

In the classical optimal control theory of stochastic systems, the value
function is a solution to the corresponding H-J-B equation whenever it has
sufficient regularity (Fleming and Rishel \cite{FR}, Krylov \cite{Krlov}).
Nevertheless when it is only known that the value function is continuous,
then, the value function is a solution to the H-J-B equation in the
viscosity sense (see Lions \cite{CIL}). From now on, we will relate the
value function of above recursive optimal singular control problem with the
following nonlinear second-order parabolic PDEs with \textquotedblleft \emph{%
gradient constraints}\textquotedblright :%
\begin{equation}
\left\{
\begin{array}{l}
\min \Big (u_{x}^{\top }\left( t,x\right) G+K,\frac{\partial }{\partial t}%
u\left( t,x\right) \\
+\min_{v\in U}\mathcal{L}\left( t,x,v\right) u\left( t,x\right) +f\left(
t,x,u\left( t,x\right) ,\nabla u\left( t,x\right) \sigma \left( t,x,v\right)
,v\right) \Big )=0, \\
u\left( T,x\right) =\Phi \left( x\right) ,\text{ }0\leq t\leq T.%
\end{array}%
\right.  \label{pde1}
\end{equation}


Before our main result, we state the following interesting result whenever
the value function defined by (\ref{value}) is sufficiently smooth.

\begin{lemma}
\label{Property}Assume that $u\left( t,x\right) \in C^{1,2}\left( \left[ 0,T%
\right] \times \mathbb{R}^{n};\mathbb{R}\right) $ is a classical solution to
H-J-B inequality \emph{(\ref{pde1}}). For any $\left( t,x\right) \in
\mathcal{D}_{t}\left( u\right) $, we have
\begin{equation}
u\left( t,x\right) \leq u\left( t,x+Gh\right) +Kh,  \label{jump}
\end{equation}%
for each $h\in \mathbb{R}_{+}^{m}.$ Furthermore, if the equality holds for
some $h\in \mathbb{R}_{+}^{m},$ then there exists $\bar{h}=\left( \bar{h}%
\right) ^{i}\in \mathbb{R}_{+}^{m}$ such that $\left( \bar{h}\right)
^{i}\leq \left( h\right) ^{i},$ $1\leq i\leq m$ and the equality also holds.
\end{lemma}

\paragraph{Proof}

From (\ref{b1}), it is easy to show%
\begin{eqnarray*}
u\left( t,x+Gh\right) -u\left( t,x\right) &=&\int_{0}^{1}u_{x}\left(
t,x+\theta Gh\right) Gh\mathrm{d}\theta \\
&\geq &-Kh.
\end{eqnarray*}%
We now proceed the second result. If the equality holds as an inequality for
some $h\in \mathbb{R}_{+}^{m},$ then we take $\bar{h}$ such that $\bar{h}%
^{i}\leq h^{i},$ $1\leq i\leq m.$ It follows that
\begin{eqnarray*}
u\left( t,x\right) -Kh &=&u\left( t,x+Gh\right) \\
&\geq &u\left( t,x+G\bar{h}\right) -K\left( h-\bar{h}\right) \\
&\geq &u\left( t,x\right) -K\bar{h}-K\left( h-\bar{h}\right) \\
&=&u\left( t,x\right) -Kh.
\end{eqnarray*}%
Thus
\begin{equation*}
u\left( t,x\right) =u\left( t,x+G\bar{h}\right) +K\bar{h}.
\end{equation*}%
The proof is completed. \hfill $\Box $

%

We now focus on the value function $u$ defined in (\ref{value}), satisfying
continuity (see Lemma \ref{LH}). Inspired by the inequality (\ref{jump}), we
may define a \emph{freedom domain} for H-J-B inequality (\ref{pde1}). The
following lemma will indicate that, for H-J-B inequality (\ref{pde1}), there
exists a region, called \textit{inaction} region, such that if the optimal
state process starts from outside this region, then apply the singular
control to bring it to the region immediately, and then keep it inside this
region from then on. In this interiors of the region, the optimal singular
control should not produce any jump.

\begin{lemma}
\label{P2}Define
\begin{equation*}
\mathcal{\tilde{D}}_{t}\left( u\right) :=\left\{ x\in \mathbb{R}^{n}:u\left(
t,x\right) <u\left( t,x+Gh\right) +Kh,\text{ }h\in \mathbb{R}_{+}^{m},\text{
}h\neq 0\right\} .
\end{equation*}%
Then the optimal state process $X^{t,x;\hat{v},\hat{\xi}}$ is continuous
whenever $\left( r,X_{r}^{t,x;\hat{v},\hat{\xi}}\right) \in \mathcal{\tilde{D%
}}_{r}\left( u\right) $. To be precise, we have%
\begin{equation*}
P\left( \Delta X_{r}^{t,x;\hat{v},\hat{\xi}}\neq 0,\text{ }X_{r}^{t,x;\hat{v}%
,\hat{\xi}}\in \mathcal{\tilde{D}}_{r}\left( u\right) \right) =0,\text{ }%
t\leq r\leq T.
\end{equation*}
\end{lemma}

\paragraph{Proof}

For any $r>t,$ we first have $\Delta X_{r}^{t,x;v,\xi }=G\cdot \Delta \xi
_{r},$ which yields that $X_{r+}^{t,x;v,\xi }=X_{r}^{t,x;v,\xi }+G\cdot
\Delta \xi _{r}.$ On the other hand, Lemma \ref{LH} implies that
\begin{equation*}
u\left( r,X_{r+}^{t,x;\hat{v},\hat{\xi}}\right) =\lim_{r^{\prime
}\downdownarrows r}u\left( r^{\prime },X_{r^{\prime }}^{t,x;\hat{v},\hat{\xi}%
}\right) .
\end{equation*}%
By DPP (Lemma \ref{DPP}), we have, for $r^{\prime }>r,$%
\begin{eqnarray}
u\left( t,x\right) &=&\mathbb{E}^{\mathscr{F}_{t}}\Bigg [\int_{t}^{r^{\prime
}}f\left( s,X_{s}^{t,x;\hat{v},\hat{\xi}},Y_{s}^{t,x;\hat{v},\hat{\xi}%
},Z_{s}^{t,x;\hat{v},\hat{\xi}},\hat{v}_{s}\right) \mathrm{d}s  \notag \\
&&+\int_{t}^{r^{\prime }}Kd\hat{\xi}_{s}+u\left( r^{\prime },X_{r^{\prime
}}^{t,x;\hat{v},\hat{\xi}}\right) \Bigg ].  \label{C1}
\end{eqnarray}%
Letting $r^{\prime }\rightarrow r,$ we immediately have%
\begin{eqnarray*}
&&\int_{t}^{r^{\prime }}f\left( s,X_{s}^{t,x;\hat{v},\hat{\xi}},Y_{s}^{t,x;%
\hat{v},\hat{\xi}},Z_{s}^{t,x;\hat{v},\hat{\xi}},\hat{v}_{s}\right) \mathrm{d%
}s \\
&\rightarrow &\int_{t}^{r}f\left( s,X_{s}^{t,x;\hat{v},\hat{\xi}},Y_{s}^{t,x;%
\hat{v},\hat{\xi}},Z_{s}^{t,x;\hat{v},\hat{\xi}},\hat{v}_{s}\right) \mathrm{d%
}s,
\end{eqnarray*}%
\begin{equation}
\int_{t}^{r^{\prime }}Kd\hat{\xi}_{s}\rightarrow \int_{t}^{r}Kd\hat{\xi}_{s}.
\label{con2}
\end{equation}%
We proceed our proof by contradiction. Suppose that
\begin{equation}
P\left( \Delta X_{r}^{t,x;v,\xi }\neq 0,\text{ }\left( r,X_{r}^{t,x;\hat{v},%
\hat{\xi}}\right) \in \mathcal{\tilde{D}}_{r}\left( u\right) \right) >0,%
\text{ }t\leq r\leq T.  \label{con3}
\end{equation}%
Then, from Lemma \ref{LH} and (\ref{C1}), we have
\begin{eqnarray*}
u\left( t,x\right) &=&\mathbb{E}^{\mathscr{F}_{t}}\Bigg [\int_{t}^{r}f\left(
s,X_{s}^{t,x;\hat{v},\hat{\xi}},Y_{s}^{t,x;\hat{v},\hat{\xi}},Z_{s}^{t,x;%
\hat{v},\hat{\xi}},\hat{v}_{s}\right) \mathrm{d}s+\int_{t}^{r}K\mathrm{d}%
\hat{\xi}_{s} \\
&&+K\cdot \Delta \hat{\xi}_{r}+u\left( r,X_{r}^{t,x;\hat{v},\hat{\xi}%
}+G\cdot \Delta \hat{\xi}_{r}\right) \Bigg ] \\
&>&\mathbb{E}^{\mathscr{F}_{t}}\Bigg [\int_{t}^{r}f\left( s,X_{s}^{t,x;\hat{v%
},\hat{\xi}},Y_{s}^{t,x;\hat{v},\hat{\xi}},Z_{s}^{t,x;\hat{v},\hat{\xi}},%
\hat{v}_{s}\right) \mathrm{d}s+\int_{t}^{r}K\mathrm{d}\hat{\xi}_{s} \\
&&+u\left( r,X_{r}^{t,x;\hat{v},\hat{\xi}}\right) \Bigg ],
\end{eqnarray*}%
which leads to a contradiction to Lemma \ref{DPP}. Note that the strict
inequity holds whenever $\left( r,X_{r}^{t,x;\hat{v},\hat{\xi}}\right) \in
\mathcal{\tilde{D}}_{r}\left( u\right) ,$ $t\leq r\leq T$, a.s. ~\hfill $%
\Box $

\begin{remark}
We say $\mathcal{\tilde{D}}_{t}\left( u\right) $ is the domain of $\xi $%
-inaction, while its complement the region of $\xi $-action. The optimal
process can be described via these two regions. For $\left(
t,X_{t}^{t,x;v,\xi }\right) =\left( t,x\right) \notin \mathcal{\tilde{D}}%
_{t}\left( u\right) ,$ the optimal control processes make $\left( t,x\right)
$ to jump instantly to the boundary of $\mathcal{\tilde{D}}_{t}\left(
u\right) .$ In fact, by Lemma \ref{Property}, if $\left( t,X_{t}^{t,x;v,\xi
}\right) \notin \mathcal{\tilde{D}}_{t}\left( u\right) $, there exists some $%
\bar{h}=\left( \bar{h}\right) ^{i}\in \mathbb{R}_{+}^{m}$ such that $\left(
\bar{h}\right) ^{i}\leq \left( h\right) ^{i},$ $1\leq i\leq m,$
\begin{equation*}
u\left( t,x\right) \leq u\left( t,x+G\bar{h}\right) +K\bar{h}.
\end{equation*}%
Hence, if $u\left( t,x\right) \in C^{1,2}\left( \left[ 0,T\right] \times
\mathbb{R}^{n};\mathbb{R}\right) $, we have $u_{x}^{\top }\left( t,x\right)
G+K\geq 0,$ which means (\ref{b1}) always holds.
\end{remark}

Our aim is to prove that the value function $u\left( t,x\right) $ introduced
by (\ref{value}), is the viscosity solution of the H-J-B variational
inequality (\ref{pde1}) of singular optimal problem. We first recall the
definition of a viscosity solution for H-J-B variational inequality (\ref%
{pde1}) from \cite{CIL}.

\begin{definition}
\label{dvis1} Let $u\left( t,x\right) \in C\left( \left[ 0,T\right] \times
\mathbb{R}^{n}\right) $ and $\left( t,x\right) \in \left[ 0,T\right] \times
\mathbb{R}^{n}.$ For every $\varphi \in C^{1,2}\left( \left[ 0,T\right]
\times \mathbb{R}^{n}\right) $

\noindent (1) for each local maximum point $\left( t_{0},x_{0}\right) $ of $%
u-\varphi $ in the interior of $\left[ 0,T\right] \times \mathbb{R}^{n},$ we
have%
\begin{equation}
\min \left( \varphi _{x}^{\top }G+K,\frac{\partial \varphi }{\partial t}%
+\min_{v\in U}\left\{ \mathcal{L}\varphi +f\left( t_{0},x_{0},\varphi
,\nabla \varphi \sigma ,v\right) \right\} \right) \geq 0  \label{vis1}
\end{equation}

at $\left( t_{0},x_{0}\right) ,$ i.e., $u$ is a subsolution.

\noindent (2) for each local minimum point $\left( t_{0},x_{0}\right) $ of $%
u-\varphi $ in the interior of $\left[ 0,T\right] \times \mathbb{R}^{n},$ we
have%
\begin{equation}
\min \left( \varphi _{x}^{\top }G+K,\frac{\partial \varphi }{\partial t}%
+\min_{v\in U}\left\{ \mathcal{L}\varphi +f\left( t_{0},x_{0},\varphi
,\nabla \varphi \sigma ,v\right) \right\} \right) \leq 0  \label{vis2}
\end{equation}

at $\left( t_{0},x_{0}\right) ,$ i.e., $u$ is a supersolution.

\noindent (3) $u\left( t,x\right) \in C\left( \left[ 0,T\right] \times
\mathbb{R}^{n}\right) $ is said to be a viscosity solution of \emph{(\ref%
{pde1})} if it is both a viscosity sub and supersolution.
\end{definition}

We have the other definition which will be useful to verify the viscosity
solutions. Below, $\mathbb{S}^{n}$ will denote the set of $n\times n$
symmetric matrices.

\begin{definition}
Let $u\left( t,x\right) \in C\left( \left[ 0,T\right] \times \mathbb{R}%
^{n}\right) $ and $\left( t,x\right) \in \left[ 0,T\right] \times \mathbb{R}%
^{n}$. We denote by $\mathcal{P}^{2,+}u\left( t,x\right) $, the
\textquotedblleft parabolic superjet\textquotedblright\ of $u$ at $\left(
t,x\right) $ the set of triples $\left( p,q,X\right) \in \mathbb{R}\times
\mathbb{R}^{n}\times \mathbb{S}^{n}$ which are such that%
\begin{eqnarray*}
u\left( s,y\right) &\leq &u\left( t,x\right) +p\left( s-t\right)
+\left\langle q,x-y\right\rangle \\
&&+\frac{1}{2}\left\langle X\left( y-x\right) ,y-x\right\rangle +o\left(
\left\vert s-t\right\vert +\left\vert y-x\right\vert ^{2}\right) .
\end{eqnarray*}%
Similarly, we denote by $\mathcal{P}^{2,-}u\left( t,x\right) ,$ the
\textquotedblleft parabolic subjet\textquotedblright\ of $u$ at $\left(
t,x\right) $ the set of triples $\left( p,q,X\right) \in \mathbb{R}\times
\mathbb{R}^{n}\times \mathbb{S}^{n}$ which are such that%
\begin{eqnarray*}
u\left( s,y\right) &\geq &u\left( t,x\right) +p\left( s-t\right)
+\left\langle q,x-y\right\rangle \\
&&+\frac{1}{2}\left\langle X\left( y-x\right) ,y-x\right\rangle +o\left(
\left\vert s-t\right\vert +\left\vert y-x\right\vert ^{2}\right) .
\end{eqnarray*}
\end{definition}

For convenience, specifically dealing with the uniqueness of viscosity
solution, we define%
\begin{eqnarray*}
\mathcal{H}\left( t,x,q,X\right) &=&\inf_{v\in U}\Big \{\frac{1}{2}\text{%
\textrm{Tr}}\left( \sigma \sigma ^{\ast }\left( t,x,v\right) X\right)
+\left\langle q,b\left( t,x,v\right) \right\rangle \\
&&+f\left( t,x,u\left( t,x\right) ,q\sigma \left( t,x,v\right) \right) \Big
\}.
\end{eqnarray*}

Note that, when the terms $\left( \frac{\partial \varphi }{\partial t}%
,\varphi _{x},\varphi _{xx}\right) $ in Definition \ref{dvis1} are taken
replace by $\left( p,q,X\right) $, we also have the other equivalent
definition of viscosity solution, that is, $\left( \frac{\partial \varphi }{%
\partial t},\varphi _{x},\varphi _{xx}\right) $ with local maximum point $%
\left( t_{0},x_{0}\right) $ of $u-\varphi $ (local minimum point $\left(
t_{0},x_{0}\right) $ of $u-\varphi $) is corresponding to $\left(
p,q,X\right) \in \mathcal{P}^{2,+}u\left( t,x\right) $ ($\left( p,q,X\right)
\in \mathcal{P}^{2,-}u\left( t,x\right) $), respectively. For more details,
see Fleming and Soner \cite{FS}, Lemma 4.1.

%

We are now in the position to assert following:

\begin{theorem}
\label{t3}The value function $u\left( \cdot ,\cdot \right) $ defined in
\emph{(\ref{value})} is a continuous viscosity solution to \emph{(\ref{pde1})%
}.
\end{theorem}

\paragraph{Proof}

From Lemma \ref{LH}, we have $u\in C\left( \left[ 0,T\right] \times \mathbb{R%
}^{n}\right) $ and clearly $u\left( T,x\right) =\Phi \left( x\right) .$\ We
first show $u$ is subsolution.To this end, we suppose that $\varphi \in
C^{1,2}\left( \left[ 0,T\right] \times \mathbb{R}^{n}\right) $, for any $%
\left( t,x\right) \in \left[ 0,T\right] \times \mathbb{R}^{n}$ such that $%
u-\varphi $ attains a global maximum at $\left( t,x\right) $ since the value
function is continuous and linear growth (see Lemma \ref{LH}). Without loss
of generality, we may also suppose that $u\left( t,x\right) =\varphi \left(
t,x\right) ,$ then%
\begin{equation}
u\left( \bar{t},\bar{x}\right) -u\left( t,x\right) \leq \varphi \left( \bar{t%
},\bar{x}\right) -\varphi \left( t,x\right) ,\text{ }\left( t,x\right) \neq
\left( \bar{t},\bar{x}\right) \in \left[ 0,T\right] \times \mathbb{R}^{n}.
\label{relation1}
\end{equation}%
Suppose that (\ref{vis1}) does not hold which means one of the follows is
valid:%
\begin{equation}
\frac{\partial \varphi }{\partial t}\left( t,x\right) +\min_{v\in U}\mathcal{%
L}\left( t,x,v\right) \varphi \left( t,x\right) +f\left( t,x,\varphi \left(
t,x\right) ,\nabla \varphi \left( t,x\right) \sigma \left( t,x,v\right)
,v\right) <0,  \label{contr11}
\end{equation}%
\begin{equation}
\varphi _{x}^{\top }\left( t,x\right) G_{t}+K<0\text{ for some }i,\text{ }%
1\leq i\leq m.  \label{contr22}
\end{equation}%
We start to deal with (\ref{contr22}) first. If (\ref{contr22}) is true for
some $i,$ then we can take $h^{i}$ small enough such that
\begin{equation*}
\varphi \left( t,x+G^{i}h^{i}\right) -\varphi \left( t,x\right) <-Kh^{i},
\end{equation*}%
where $G^{i}$ denotes the $i$th column of $G$ ($n\times m$ matrix). Hence, (%
\ref{relation1}) can imply that
\begin{equation*}
u\left( t,x+G^{i}h^{i}\right) -u\left( t,x\right) <-Kh^{i}.
\end{equation*}%
If we set $h=\left( 0,\ldots ,h^{i},\ldots 0\right) $, then we have $u\left(
t,x\right) >u\left( t,x+Gh\right) +Kh,$ which is obviously contradiction to (%
\ref{jump}) in Lemma \ref{Property}. Now we have obtained
\begin{equation}
\varphi _{x}^{\top }\left( t,x\right) G+K\geq 0,\text{ }\left( t,x\right)
\in \left[ 0,T\right] \times \mathbb{R}^{n}.  \label{rr1}
\end{equation}%
Based on this result, we are able to construct the second contradiction for (%
\ref{contr11}). Thus by virtue of DPP (see Lemma \ref{DPP}), we have%
\begin{equation*}
\varphi \left( t,x\right) =u\left( t,x\right) =ess\inf_{\left( v,\xi \right)
\in \mathcal{U}}\mathcal{G}_{t,t+\delta }^{t,x;v,\xi }\left[ u\left(
t+\delta ,X_{t+\delta }^{t,x;v,\xi }\right) \right] ,\text{ }0\leq \delta
\leq T-t.
\end{equation*}%
To abbreviate the notations, we put$,$%
\begin{eqnarray*}
F\left( s,x,y,z,v\right) &=&\frac{\partial }{\partial s}\varphi \left(
s,x\right) +\mathcal{L}\left( t,x,v\right) \varphi \left( t,x\right) \\
&&+f\left( s,x,y+\varphi \left( s,x\right) ,z+\nabla \varphi \left(
s,x\right) \sigma \left( s,x,v\right) ,v\right) ,
\end{eqnarray*}%
where $\left( s,x,y,z,v\right) \in \left[ 0,T\right] \times \mathbb{R}%
^{n}\times \mathbb{R\times R}^{d}\times U.$ We now introduce the following
BSDE on interval $\left[ t,t+\delta \right] $, $0\leq \delta \leq T-t,$%
\begin{equation}
Y_{s}^{1,v,\xi }=\int_{s}^{t+\delta }F\left( s,X_{s}^{t,x;v,\xi
},Y_{s}^{1,v,\xi },Z_{s}^{1,v,\xi },v_{s}\right) \mathrm{d}%
s-\int_{s}^{t+\delta }Z_{s}^{1,v,\xi }\mathrm{d}W_{s}+\int_{s}^{t+\delta }K%
\mathrm{d}\xi _{r},  \label{BV1}
\end{equation}%
where $\left( v,\xi \right) \in \mathcal{U}$. Clearly $Y_{t+\delta
}^{1,v,\xi }=0.$ We check that $F\left( s,X_{s}^{t,x;v,\xi
},y,z,v_{s}\right) $ satisfies the assumption (A2), thus we claim that Eq. (%
\ref{BV1}) admits a unique adapted solution. In order to characterize the
process $Y_{s}^{1,v,\xi }$, we will study a new BSDE as follows:%
\begin{eqnarray}
Y_{s}^{v,\xi } &=&\varphi \left( t+\delta ,X_{t+\delta }^{t,x;v,\xi }\right)
+\int_{s}^{t+\delta }f\left( r,X_{r}^{t,x;v,\xi },Y_{r}^{v,\xi
},Z_{r}^{v,\xi },v_{r}\right) \mathrm{d}r  \notag \\
&&-\int_{s}^{t+\delta }Z_{r}^{v,\xi }\mathrm{d}W_{s}+\int_{s}^{t+\delta }K%
\mathrm{d}\xi _{r}.  \label{BV3}
\end{eqnarray}%
Applying It\^{o}'s formula to $\varphi \left( t,x\right) $ on interval $%
\left[ t,t+\delta \right] ,$ we have
\begin{eqnarray}
\varphi \left( s,X_{s}^{t,x;v,\xi }\right) &=&\varphi \left( t+\delta
,X_{t+\delta }^{t,x;v,\xi }\right)  \notag \\
&&-\int_{s}^{t+\delta }\left[ \frac{\partial }{\partial s}\varphi \left(
r,x\right) +\mathcal{L}\left( r,X_{r}^{t,x;v,\xi },v_{r}\right) \varphi
\left( r,X_{r}^{t,x;v,\xi }\right) \right] \mathrm{d}r  \notag \\
&&-\int_{s}^{t+\delta }\nabla \varphi ^{\top }\left( r,X_{r}^{t,x;v,\xi
}\right) G\mathrm{d}\xi _{r}  \notag \\
&&-\int_{s}^{t+\delta }\nabla \varphi \left( r,X_{r}^{t,x;v,\xi }\right)
\sigma \left( r,X_{r}^{t,x;v,\xi },v_{r}\right) \mathrm{d}W_{r}  \notag \\
&&-\sum_{s\leq r\leq t+\delta }\left[ \varphi \left( r,X_{r+}^{t,x;v,\xi
}\right) -\varphi \left( r,X_{r}^{t,x;v,\xi }\right) \right.  \notag \\
&&\left. -\nabla \varphi ^{\top }\left( r,X_{r}^{t,x;v,\xi }\right) \Delta
X_{r}^{t,x;v,\xi }\right] .  \label{CB}
\end{eqnarray}%
Since
\begin{eqnarray*}
&&-\int_{s}^{t+\delta }\nabla \varphi ^{\top }\left( r,X_{r}^{t,x;v,\xi
}\right) G\mathrm{d}\xi _{r}+\sum_{s\leq r\leq t+\delta }\nabla \varphi
^{\top }\left( r,X_{r}^{t,x;v,\xi }\right) \Delta X_{r}^{t,x;v,\xi } \\
&=&-\int_{s}^{t+\delta }\nabla \varphi ^{\top }\left( r,X_{r}^{t,x;v,\xi
}\right) G\mathrm{d}\xi _{r}^{c} \\
&\leq &\int_{s}^{t+\delta }K\mathrm{d}\xi _{r}^{c}
\end{eqnarray*}%
and
\begin{equation*}
\varphi \left( r,X_{r+}^{t,x;v,\xi }\right) -\varphi \left(
r,X_{r}^{t,x;v,\xi }\right) \geq -K\Delta \xi _{r},
\end{equation*}%
we have%
\begin{eqnarray*}
\varphi \left( s,X_{s}^{t,x;v,\xi }\right) &\leq &\varphi \left( t+\delta
,X_{t+\delta }^{t,x;v,\xi }\right) \\
&&-\int_{s}^{t+\delta }\frac{\partial }{\partial s}\varphi \left( r,x\right)
+\mathcal{L}\left( r,X_{r}^{t,x;v,\xi },v_{r}\right) \varphi \left(
r,X_{r}^{t,x;v,\xi }\right) \\
&&-\int_{s}^{t+\delta }\nabla \varphi \left( r,X_{r}^{t,x;v,\xi }\right)
\sigma \left( r,X_{r}^{t,x;v,\xi },v_{r}\right) \mathrm{d}W_{r} \\
&&+\int_{s}^{t+\delta }K\mathrm{d}\xi _{r},
\end{eqnarray*}%
only when $\xi _{r}\equiv 0$ above takes equality. For simplicity, putting $%
\xi \equiv 0,$ from (\ref{BV3}) and (\ref{CB}), we can get
\begin{eqnarray*}
Y_{s}^{v,0}-\varphi \left( s,X_{s}^{t,x;v,0}\right) &=&\int_{s}^{t+\delta }
\left[ f\left( r,X_{r}^{t,x;v,0},Y_{r}^{v,0},Z_{r}^{v,0},v_{r}\right) \right.
\\
&&\left. +\frac{\partial }{\partial s}\varphi \left( r,x\right) +\mathcal{L}%
\left( r,X_{r}^{t,x;v,0},v_{r}\right) \varphi \left(
r,X_{r}^{t,x;v,0}\right) \right] \mathrm{d}r \\
&&-\int_{s}^{t+\delta }\left[ Z_{r}^{v,\xi }-\nabla \varphi \left(
r,X_{r}^{t,x;v,0}\right) \sigma \left( r,X_{r}^{t,x;v,0},v_{r}\right) \right]
\mathrm{d}W_{r},
\end{eqnarray*}%
which means that
\begin{eqnarray*}
Y_{s}^{1,v,0} &=&Y_{s}^{v,0}-\varphi \left( s,X_{s}^{t,x;v,0}\right) \\
&=&\mathcal{G}_{s,t+\delta }^{t,x;v,0}\left[ \varphi \left( t+\delta
,X_{t+\delta }^{t,x;v,0}\right) \right] -\varphi \left(
s,X_{s}^{t,x;v,0}\right) ,t\leq s\leq t+\delta .
\end{eqnarray*}%
Note that in the coefficient $F$ in Eq. (\ref{BV1}), it contains the
stochastic term $X_{s}^{t,x;v,\xi }$. We want to take replace it by $x\in
\mathbb{R}^{n}.$ So we consider the following BSDE also on interval $\left[
t,t+\delta \right] $ with terminal value $0,$%
\begin{equation}
Y_{s}^{2,v}=\int_{s}^{t+\delta }F\left(
s,x,Y_{s}^{2,v},Z_{s}^{2,v},v_{s}\right) \mathrm{d}s-\int_{s}^{t+\delta
}Z_{s}^{2,v}\mathrm{d}W_{s},  \label{BV2}
\end{equation}%
where $v\in \mathcal{U}_{1}$. From Lemma \ref{state4} (see in Appendix), we
have
\begin{equation}
\left\vert Y_{s}^{1,v,0}-Y_{s}^{2,v}\right\vert \leq C\delta ^{\frac{3}{2}},%
\text{ }t\leq s\leq t+\delta .  \label{e1}
\end{equation}%
By virtue of DPP (see Lemma \ref{DPP}), we have%
\begin{equation*}
\varphi \left( t,x\right) =u\left( t,x\right) =ess\inf_{\left( v,\xi \right)
\in \mathcal{U}}\mathcal{G}_{t,t+\delta }^{t,x;v,\xi }\left[ u\left(
t+\delta ,X_{t+\delta }^{t,x;v,\xi }\right) \right] ,\text{ }0\leq \delta
\leq T-t.
\end{equation*}%
From $u\leq \varphi $ and the monotonicity property of $G_{t,t+\delta
}^{t,x;v,\xi }\left[ \cdot \right] $, it follows%
\begin{eqnarray*}
\varphi \left( t,x\right) &=&u\left( t,x\right) \\
&=&ess\inf_{\left( v,\xi \right) \in \mathcal{U}}\mathcal{G}_{t,t+\delta
}^{t,x;v,\xi }\left[ u\left( t+\delta ,X_{t+\delta }^{t,x;v,\xi }\right) %
\right] \\
&\leq &ess\inf_{\left( v,\xi \right) \in \mathcal{U}}\mathcal{G}_{t,t+\delta
}^{t,x;v,\xi }\left[ \varphi \left( t+\delta ,X_{t+\delta }^{t,x;v,\xi
}\right) \right] \\
&\leq &ess\inf_{v\in \mathcal{U}_{1}}\mathcal{G}_{t,t+\delta }^{t,x;v,0}%
\left[ \varphi \left( t+\delta ,X_{t+\delta }^{t,x;v,0}\right) \right] ,
\end{eqnarray*}%
which yields
\begin{equation}
ess\inf_{v\in \mathcal{U}_{1}}Y_{t}^{1,v,0}\geq 0.  \label{e2}
\end{equation}%
From (\ref{e1}), we get%
\begin{equation}
Y_{t}^{2,v}\geq Y_{t}^{1,v,0}-C\delta ^{\frac{3}{2}}.  \label{e3}
\end{equation}%
Combining to (\ref{e2}) and (\ref{e3}), we must have%
\begin{equation}
Y_{t}^{2,v}\geq -C\delta ^{\frac{3}{2}}.  \label{est1}
\end{equation}%
Observe that in Eq. (\ref{BV2}), $F$ depends on control $v$. Now, eliminate
the control term via introducing the following BSDE:%
\begin{equation}
\mathcal{Y}_{s}=\int_{s}^{t+\delta }F_{0}\left( s,x,\mathcal{Y}_{s},\mathcal{%
Z}_{s}\right) \mathrm{d}s-\int_{s}^{t+\delta }\mathcal{Z}_{s}\mathrm{d}W_{s},
\label{BV4}
\end{equation}%
where
\begin{equation}
F_{0}\left( s,x,y,z\right) =\inf_{v\in U}F\left( s,x,y,z,v\right) .
\label{F0}
\end{equation}%
Clearly,
\begin{equation}
F\left( s,x,y,z,v\right) \geq F_{0}\left( s,x,y,z\right) ,  \label{comf}
\end{equation}%
where $\left( s,x,y,z,v\right) \in \left[ 0,T\right] \times \mathbb{R}%
^{n}\times \mathbb{R\times R}^{d}\times U.$ Applying the comparison theorem
(Lemma \ref{comp}) to Eq. (\ref{BV2}) and (\ref{BV4}), we have%
\begin{equation}
ess\inf_{v\in \mathcal{U}_{1}}Y_{t}^{2,v}\geq \mathcal{Y}_{t},\text{ a.s..}
\label{cx1}
\end{equation}%
Indeed, from (\ref{comf}) and Lemma \ref{comp}, it follows that $%
Y_{t}^{2,v}\geq \mathcal{Y}_{t},$ a.s. for any $v\in \mathcal{U}_{1}.$
Taking $essinf$, we get (\ref{cx1}) immediately. Besides, there exists a
measurable function $v^{\prime }:\left[ t,T\right] \times \mathbb{R}%
^{n}\times \mathbb{R\times R}^{d}\rightarrow U$ such that
\begin{equation*}
F_{0}\left( s,x,y,z\right) =F\left( s,x,y,z,v^{\prime }\left( s,x,y,z\right)
\right) ,\text{ }\left( s,x,y,z\right) \in \left[ 0,T\right] \times \mathbb{R%
}^{n}\times \mathbb{R\times R}^{d}.
\end{equation*}%
Then taking $v_{s}^{\prime }:=v^{\prime }\left( s,x,\mathcal{Y}_{s},\mathcal{%
Z}_{s}\right) ,$ $t\leq s\leq t+\delta ,$ we have $v_{\cdot }^{\prime }\in
\mathcal{U}_{1}$ and
\begin{equation*}
F_{0}\left( s,x,\mathcal{Y}_{s},\mathcal{Z}_{s}\right) =F\left( s,x,\mathcal{%
Y}_{s},\mathcal{Z}_{s},v_{s}^{\prime }\right) ,\text{ }t\leq s\leq t+\delta .
\end{equation*}%
According to the uniqueness solution of BSDE (\ref{BV2}) and (\ref{BV4}), it
follows that $\left( \mathcal{Y}_{s},\mathcal{Z}_{s}\right) =\left(
Y_{s}^{2,v^{\prime }},Z_{s}^{2,v^{\prime }}\right) .$ Particularly,
\begin{equation}
\mathcal{Y}_{t}=Y_{t}^{2,v^{\prime }},\text{ a.s..}  \label{cx2}
\end{equation}%
From (\ref{cx1}) and (\ref{cx2}), it follows that $ess\inf_{v\times \mathcal{%
U}_{1}}Y_{t}^{2,v}=\mathcal{Y}_{t}.$ Then (\ref{BV4}) becomes
\begin{equation}
\mathcal{Y}_{s}^{0}=\int_{s}^{t+\delta }F_{0}\left( s,x,\mathcal{Y}_{s}^{0},%
\mathcal{Z}_{s}^{0}\right) \mathrm{d}s-\int_{s}^{t+\delta }\mathcal{Z}%
_{s}^{0}\mathrm{d}W_{s}.  \label{BV5}
\end{equation}%
Moreover, from the definition of $F_{0},$ it follows that $F_{0}\left(
s,x,y,z\right) $ is deterministic function. Note that the terminal value is $%
0,$ we immediately get the solution to Eq. (\ref{BV5}) is a triple $\left(
\mathcal{Y}_{s}^{0},\mathcal{Z}_{s}^{0}\right) =\left( \mathcal{\bar{Y}}%
_{s}^{0},0\right) ,$ where $\mathcal{\bar{Y}}_{s}^{0}$ satisfies the
following ODE:%
\begin{equation}
\mathcal{\bar{Y}}_{t}^{0}=\int_{t}^{t+\delta }F_{0}\left( s,x,\mathcal{\bar{Y%
}}_{s}^{0},0\right) \mathrm{d}s.  \label{ode1}
\end{equation}%
On the other hand, (\ref{est1}) implies that $\int_{t}^{t+\delta
}F_{0}\left( s,x,\mathcal{\bar{Y}}_{s}^{0},0\right) \mathrm{d}s\geq -C\delta
^{\frac{3}{2}}.$ By classical estimate$,$ we obtain $F_{0}\left(
s,x,0,0\right) \geq 0,$ that is
\begin{eqnarray*}
0 &\leq &F_{0}\left( s,x,0,0\right) \\
&=&\inf_{v\in U}F\left( s,x,0,0,v\right) \\
&=&\frac{\partial \varphi }{\partial t}\left( t,x\right) \\
&&+\min_{v\in U}\mathcal{L}\left( t,x,v\right) \varphi \left( t,x\right)
+f\left( t,x,\varphi \left( t,x\right) ,\nabla \varphi \left( t,x\right)
\sigma \left( t,x,v\right) ,v\right) ,
\end{eqnarray*}%
which leads contradiction to (\ref{contr11}).

Next we show that $u$ is meanwhile a supersolution of (\ref{pde1}). Suppose
that $\varphi \in C^{1,2}\left( \left[ 0,T\right] \times \mathbb{R}%
^{n}\right) $, for any $\left( t,x\right) \in \left[ 0,T\right] \times
\mathbb{R}^{n}$ such that $u-\varphi $ attains a global minimum at $\left(
t,x\right) $. Without loss of generality we may also suppose that $u\left(
t,x\right) =\varphi \left( t,x\right) ,$ then%
\begin{equation}
u\left( \bar{t},\bar{x}\right) -u\left( t,x\right) \geq \varphi \left( \bar{t%
},\bar{x}\right) -\varphi \left( t,x\right) ,\text{ }\left( t,x\right) \neq
\left( \bar{t},\bar{x}\right) \in \left[ 0,T\right] \times \mathbb{R}^{n}.
\label{relation2}
\end{equation}%
If (\ref{vis2}) fails, the following inequity holds true simultaneously,
\begin{equation}
\frac{\partial \varphi }{\partial t}\left( t,x\right) +\min_{v\in U}\mathcal{%
L}\left( t,x,v\right) \varphi \left( t,x\right) +f\left( t,x,\varphi \left(
t,x\right) ,\nabla \varphi \left( t,x\right) \sigma \left( t,x,v\right)
,v\right) >0,  \label{sup1}
\end{equation}%
\begin{equation}
\varphi _{x}^{\top }\left( t,x\right) G_{t}+K>0,\text{ for }i,\text{ }1\leq
i\leq m.  \label{sup2}
\end{equation}%
For small $h\in \mathbb{R}_{+}^{m}$, $h\neq 0,$ $\varphi \left(
t,x+Gh\right) -\varphi \left( t,x\right) \geq -Kh.$ Thus, by (\ref{relation2}%
), $u\left( t,x+Gh\right) -u\left( t,x\right) \geq -Kh,$ which means $\left(
r,x\right) \in \mathcal{D}\left( u\right) ,$ and by Lemma \ref{P2}, we know
that $\hat{X}^{t,x;\hat{v},\hat{\xi}}$ is continuous for optimal controls $%
\left( \hat{v},\hat{\xi}\right) \in \mathcal{U}$. Hence, applying It\^{o}'s
formula to $\varphi \left( t,x\right) $, we have%
\begin{eqnarray}
\varphi \left( s,X_{s}^{t,x;\hat{v},\hat{\xi}}\right) &=&\varphi \left(
t+\delta ,X_{t+\delta }^{t,x;\hat{v},\hat{\xi}}\right)  \notag \\
&&-\int_{s}^{t+\delta }\left[ \frac{\partial }{\partial s}\varphi \left(
r,x\right) +\mathcal{L}\left( r,X_{r}^{t,x;\hat{v},\hat{\xi}},\hat{v}%
_{r}\right) \varphi \left( r,X_{r}^{t,x;\hat{v},\hat{\xi}}\right) \right]
\mathrm{d}r  \notag \\
&&-\int_{s}^{t+\delta }\nabla \varphi ^{\top }\left( r,X_{r}^{t,x;\hat{v},%
\hat{\xi}}\right) G\mathrm{d}\hat{\xi}_{r}  \notag \\
&&-\int_{s}^{t+\delta }\nabla \varphi \left( r,X_{r}^{t,x;\hat{v},\hat{\xi}%
}\right) \sigma \left( r,X_{r}^{t,x;\hat{v},\hat{\xi}},\hat{v}_{r}\right)
\mathrm{d}W_{r}.  \label{superito}
\end{eqnarray}%
Repeating the method above, we introduce the following two BSDEs:%
\begin{eqnarray}
Y_{s}^{\hat{v},\hat{\xi}} &=&\varphi \left( t+\delta ,X_{t+\delta }^{t,x;%
\hat{v},\hat{\xi}}\right)  \notag \\
&&+\int_{s}^{t+\delta }f\left( r,X_{r}^{t,x;\hat{v},\hat{\xi}},Y_{r}^{\hat{v}%
,\hat{\xi}},Z_{r}^{\hat{v},\hat{\xi}},\hat{v}_{r}\right) \mathrm{d}r  \notag
\\
&&-\int_{s}^{t+\delta }Z_{r}^{\hat{v},\hat{\xi}}\mathrm{d}%
W_{s}+\int_{s}^{t+\delta }K\mathrm{d}\hat{\xi}_{r}  \label{super1}
\end{eqnarray}%
and%
\begin{eqnarray}
Y_{s}^{1,\hat{v},\hat{\xi}} &=&\int_{s}^{t+\delta }F\left( s,X_{s}^{t,x;\hat{%
v},\hat{\xi}},Y_{s}^{1,\hat{v},\hat{\xi}},Z_{s}^{1,\hat{v},\hat{\xi}},\hat{v}%
_{s}\right) \mathrm{d}s-\int_{s}^{t+\delta }Z_{s}^{1,\hat{v},\hat{\xi}}%
\mathrm{d}W_{s}  \notag \\
&&+\int_{s}^{t+\delta }\left( K+\nabla \varphi ^{\top }\left( r,X_{r}^{t,x;%
\hat{v},\hat{\xi}}\right) G_{r}\right) \mathrm{d}\hat{\xi}_{r}.
\label{super2}
\end{eqnarray}%
Clearly,
\begin{equation}
Y_{s}^{1,\hat{v},\hat{\xi}}=Y_{s}^{\hat{v},\hat{\xi}}-\varphi \left(
s,X_{s}^{t,x;\hat{v},\hat{\xi}}\right) .  \label{dpp2}
\end{equation}%
Consider the following BSDE:
\begin{eqnarray}
Y_{s}^{2,\hat{v},\hat{\xi}} &=&\int_{s}^{t+\delta }F\left( s,x,Y_{s}^{2,\hat{%
v},\hat{\xi}},Z_{s}^{2,\hat{v},\hat{\xi}},\hat{v}_{s}\right) \mathrm{d}%
s-\int_{s}^{t+\delta }Z_{s}^{2,\hat{v},\hat{\xi}}\mathrm{d}W_{s}  \notag \\
&&+\int_{s}^{t+\delta }\left( K+\nabla \varphi ^{\top }\left( r,X_{r}^{t,x;%
\hat{v},\hat{\xi}}\right) G_{r}\right) \mathrm{d}\hat{\xi}_{r}.
\label{super3}
\end{eqnarray}%
From Lemma \ref{adcontrol} (see in Appendix), we have
\begin{equation}
\left\vert Y_{s}^{1,\hat{v},\hat{\xi}}-Y_{s}^{2,\hat{v},\hat{\xi}%
}\right\vert \leq C\delta ^{\frac{3}{2}},\text{ }t\leq s\leq t+\delta ,
\label{super4}
\end{equation}%
which implies that
\begin{equation}
0\leq Y_{t}^{2,\hat{v},\hat{\xi}}-C\delta ^{\frac{3}{2}}\leq Y_{t}^{1,\hat{v}%
,\hat{\xi}}.  \label{superkey}
\end{equation}%
Now we give the following BSDE:
\begin{eqnarray}
\mathcal{Y}_{s}^{0,\hat{\xi}} &=&\int_{s}^{t+\delta }F_{0}\left( s,x,%
\mathcal{Y}_{s}^{0,\hat{\xi}},\mathcal{Z}_{s}^{0,\hat{\xi}}\right) \mathrm{d}%
s-\int_{s}^{t+\delta }\mathcal{Z}_{s}^{0,\hat{\xi}}\mathrm{d}W_{s}  \notag \\
&&+\int_{s}^{t+\delta }\left( K+\nabla \varphi ^{\top }\left( r,X_{r}^{t,x;%
\hat{v},\hat{\xi}}\right) G\right) \mathrm{d}\hat{\xi}_{r},  \label{super5}
\end{eqnarray}%
where $F_{0}\left( s,x,y,z\right) $ is defined in (\ref{F0}). By comprison
theorem (Lemma \ref{comp}), we have $\mathcal{Y}_{s}^{0,\hat{\xi}}\leq
Y_{s}^{2,\hat{v},\hat{\xi}},$ a.s. $s\leq t\leq T.$ Since $\hat{v}\in
\mathcal{U}_{1}$, we can find $\bar{v}\left( s,x,y,z\right) \in \arg
\min_{v\in U}F\left( s,x,y,z,v\right) ,$ such that $\hat{v}_{s}=\bar{v}%
\left( s,X_{r}^{t,x;\hat{v},\hat{\xi}},Y_{r}^{t,x;\hat{v},\hat{\xi}%
},Z_{r}^{t,x;\hat{v},\hat{\xi}}\right) $ and $F_{0}\left( s,x,\mathcal{Y}%
_{s}^{0,\hat{\xi}},\mathcal{Z}_{s}^{0,\hat{\xi}}\right) =F\left( s,x,%
\mathcal{Y}_{s}^{0,\hat{\xi}},\mathcal{Z}_{s}^{0,\hat{\xi}},\hat{v}%
_{s}\right) ,$ $s\in \left[ t,T\right] .$ Therefore, $\mathcal{Y}_{s}^{0,%
\hat{\xi}}=Y_{s}^{2,\hat{v},\hat{\xi}},$ a.s. $s\leq t\leq T.$ Particularly,
$\mathcal{Y}_{t}^{0,\hat{\xi}}=Y_{t}^{2,\hat{v},\hat{\xi}},$ a.s. From (\ref%
{superkey}), we have $0<\mathcal{Y}_{t}^{0,\hat{\xi}}-C\delta ^{\frac{3}{2}%
}\leq Y_{t}^{1,\hat{v},\hat{\xi}}.$ Observing (\ref{sup2}), we have $%
\mathcal{Y}_{t}^{0,0}\leq \mathcal{Y}_{t}^{0,\hat{\xi}},$where $\mathcal{Y}%
_{t}^{0,0}$ satisfies
\begin{equation}
\mathcal{Y}_{s}^{0,0}=\int_{s}^{t+\delta }F_{0}\left( s,x,\mathcal{Y}%
_{s}^{0,0},\mathcal{Z}_{s}^{0,0}\right) \mathrm{d}s-\int_{s}^{t+\delta }%
\mathcal{Z}_{s}^{0,0}\mathrm{d}W_{s}.  \label{super6}
\end{equation}%
Eq.(\ref{super6}) admits a unique adapted solution $\left( \mathcal{Y}_{s}^{0,0},%
\mathcal{Z}_{s}^{0,0}\right) =\left( \mathcal{\tilde{Y}}_{s}^{0},0\right) ,$
where $\mathcal{\tilde{Y}}_{s}^{0}$ satisfies $\mathcal{\tilde{Y}}%
_{t}^{0}=\int_{t}^{t+\delta }F_{0}\left( s,x,\mathcal{\tilde{Y}}%
_{s}^{0},0\right) \mathrm{d}s.$ Hence, we have $0<\mathcal{\tilde{Y}}%
_{s}^{0}-C\delta ^{\frac{3}{2}}\leq Y_{t}^{1,\hat{v},\hat{\xi}}.$ From (\ref%
{sup1}), we have $Y_{t}^{1,\hat{v},\hat{\xi}}>0,$ which contradicts the
dynamic programming principle. Indeed,
\begin{eqnarray*}
\varphi \left( t,x\right) &=&u\left( t,x\right) \\
&=&\mathcal{G}_{t,t+\delta }^{t,x;\hat{v},\hat{\xi}}\left[ u\left( t+\delta
,X_{t+\delta }^{t,x;\hat{v},\hat{\xi}}\right) \right] \\
&\geq &\mathcal{G}_{t,t+\delta }^{t,x;\hat{v},\hat{\xi}}\left[ \varphi
\left( t+\delta ,X_{t+\delta }^{t,x;\hat{v},\hat{\xi}}\right) \right] ,
\end{eqnarray*}%
which means that $Y_{t}^{1,\hat{v},\hat{\xi}}\leq 0.$ The proof is complete.
~\hfill $\Box $

In order to establish a uniqueness result for viscosity solution of (\ref%
{pde1}), we adapt some techniques and methods from \cite{FS}. We have to
mention that there is another approach developed by Barles, Buckdahn, and
Pardoux \cite{BBP}. The value function can be considered in given class of
continuous functions satisfying
\begin{equation*}
\lim_{\left\vert x\right\vert \rightarrow \infty }\left\vert u\left(
t,x\right) \right\vert \exp \left\{ -A\left[ \log \left( \left\vert
x\right\vert \right) \right] ^{2}\right\} =0,
\end{equation*}%
uniformly for $t\in \left[ 0,T\right] ,$ for some $A>0.$ These techniques
and methods can also be found in \cite{BJ, WY, CF} for the uniqueness for
viscosity solutions of recursive control of the obstacle constraint problem
and Hamilton-Jacobi-Bellman-Isaacs equations related to stochastic
differential games, respectively.

However, as you may have observed, in our H-J-B inequality, there appears a
term $u_{x}^{\top }\left( t,x\right) G+K$, which is a row vector (in
contrast to \cite{WY}, scalar there). It is worth to pointing out that the
smooth supersolution built in \cite{BBP}, namely
\begin{equation*}
\chi \left( t,x\right) =\exp \left[ \left( \check{C}\left( T-t\right)
+A\right) \psi \left( x\right) \right] ,
\end{equation*}%
whilst
\begin{equation*}
\psi \left( x\right) =\left[ \log \left( \left\vert x\right\vert
^{2}+1\right) ^{\frac{1}{2}}+1\right] ^{2},
\end{equation*}%
where $\check{C}$ and $A$ are positive constants has some properties. After
some tedious computations, we have%
\begin{equation*}
\chi _{x}\left( t,x\right) =2\chi \left( t,x\right) \left( \check{C}\left(
T-t\right) +A\right) \cdot \left( \left\vert x\right\vert ^{2}+1\right) ^{-%
\frac{3}{2}}x.
\end{equation*}%
It follows
\begin{equation*}
\chi _{x}^{T}\left( t,x\right) G=2\chi \left( t,x\right) \left( \check{C}%
\left( T-t\right) +A\right) \cdot \left( \left\vert x\right\vert
^{2}+1\right) ^{-\frac{3}{2}}x^{\top }\cdot G.
\end{equation*}%
Obviously, for any $x\in \mathbb{R}^{n},$ it is impossible to get $x^{\top
}\cdot G\geq 0,$ unless selecting $x$ such that $x^{\top }\cdot G\geq 0.$

Consider the following
\begin{equation}
\left\{
\begin{array}{l}
\min \Big (w_{x}^{\top }\left( t,x\right) G,\text{ } \\
\frac{\partial w}{\partial t}\left( t,x\right) +\min_{v\in U}\left\{
\mathcal{L}\left( t,x,v\right) w\left( t,x\right) +L\left\vert w\right\vert
+L\left\vert \nabla w\sigma \left( t,x,v\right) \right\vert \right\} \Big )%
=0, \\
w\left( T,x\right) =0,\text{ }0\leq t\leq T.%
\end{array}%
\right.  \label{bsuper}
\end{equation}%
On the one hand, repeating the method in \cite{BBP}, one can show
\begin{equation*}
\min_{v\in U}\left\{ \mathcal{L}\left( t,x,v\right) \chi \left( t,x\right)
+C\left\vert \chi \right\vert +C\left\vert \nabla \chi \sigma \left(
t,x,v\right) \right\vert \right\} \leq 0,
\end{equation*}%
where $C$ is the Lipschitz constant of $f.$ On the other hand, when choosing
$x\in \mathcal{A}:=\left\{ x:x^{\top }\cdot G\geq 0\right\} ,$ $\chi $ is
indeed a supersolution of (\ref{bsuper}) according to (ii) of Definition \ref%
{dvis1}. However, following the idea in \cite{BBP}, whenever considering the
difference of $u_{1}-u_{2}$ where $u_{1}$ ($u_{2}$) is a subsolution
(supersolution) of (\ref{pde1}). The domain of $u_{1}-u_{2}$ must be
restricted in $\mathcal{A}$. Hence, we borrow the idea from Fleming, Soner
\cite{FS} and Haussmann, Suo \cite{HS2} to handle the uniqueness.

For our aim, let us define $\Sigma =\left[ 0,T\right] \times \mathbb{R}^{n}$
and $C\left( \Sigma ;\mathbb{R}\right) $ denote the collection of
real-valued continuous functions defined on $\Sigma ,$ value in $\mathbb{R}.$
Now introduce the function space%
\begin{equation*}
\begin{array}{lll}
\mathcal{C}\left( \Sigma \right) & := & \Big \{u:u\in C\left( \Sigma ;%
\mathbb{R}\right) ,\text{ with }u\text{ bounded, and} \\
&  & \left\vert u\left( t,x\right) -u\left( t,y\right) \leq C\left\vert
x-y\right\vert \right\vert \text{ for some }C>0\Big \}.%
\end{array}%
\end{equation*}%
We now give the uniqueness result for (\ref{pde1}). This result is attained
under more restrictive assumptions that the existence one. We add the following two additional assumptions:

\begin{description}
\item[(A4)] Assume that $b$, $\sigma ,$ $\Phi ,$ $f\left( t,\cdot
,y,z,v\right) $ are uniformly continuous, uniformly with respect to $\left(
t,y,z,v\right) ,$ and bounded.

\item[(A5)] Assume that
\begin{equation*}
\left\vert f\left( t,x,r,p,v\right) -f\left( t,x,y,p,v\right) \right\vert
\leq \varpi _{R}\left( \left\vert x-y\right\vert \left( 1+\left\vert
p\right\vert \right) \right)
\end{equation*}%
where $\varpi _{R}\left( s\right) \rightarrow 0,$ when $s\rightarrow 0^{+},$
for all $t\in \left[ 0,T\right] ,$ $\left\vert x\right\vert ,$ $\left\vert
y\right\vert \leq R,$ $\left\vert r\right\vert \leq R,$ $\forall R>0.$
\end{description}

We have put somewhat strong assumptions, namely, $b,$ $\sigma ,$ $f$ are
bounded. These conditions may be removed by modifying the idea by Ishii \cite%
{Is}.

\begin{theorem}
\label{t4}Assume that \emph{(A1)-(A5)} are in force. Then there exists at
most one viscosity solution of H-J-B inequality \emph{(\ref{pde1})} in the
class of bounded and continuous functions.
\end{theorem}

\paragraph{Proof}

Under the assumptions (A1)-(A4), by Proposition 2.5 in \cite{BBP} and Lemma %
\ref{LH}, one can show that the value function $u$ is uniformly continuous
and bounded. Thus, $u$ $\left( \text{defined by }(\ref{value}\right) )\in
\mathcal{C}\left( \Sigma \right) .$

Let $u_{1}$ and $u_{2}$ be two viscosity solution to (\ref{pde1}). We will
show that $u_{1}\leq u_{2}$ on $\Sigma .$ For any given $\epsilon
:0<\epsilon <1,$ define $\Sigma _{\epsilon }:=\left] \epsilon ,T\right]
\times \mathbb{R}^{n}.$ For any $\left( t,x\right) \in \Sigma _{\epsilon },$%
\begin{equation*}
u_{1}^{\epsilon }\left( t,x\right) =\left( 1-\epsilon \right) u_{1}\left(
t,x\right) -\frac{\epsilon }{1-\epsilon }.
\end{equation*}%
It is easy to check that $\frac{1}{1-\epsilon }u_{1}^{\epsilon }$ is a
viscosity subsolution to (\ref{pde1}). Given any $\alpha ,$ $\beta >0,$ $%
0<\epsilon <1,$ define an auxiliary function%
\begin{equation*}
\Psi \left( t,x,y\right) =u_{1}^{\epsilon }\left( t,x\right) -u_{2}\left(
t,y\right) -\frac{1}{2\alpha }\left\vert x-y\right\vert ^{2}+\beta \left(
t-T\right) .
\end{equation*}
Clearly, $\Psi $\ is bounded above on $\Sigma _{\epsilon }.$ We can find $%
\left( t_{\alpha },x_{\alpha },y_{\alpha }\right) \in \Sigma _{\epsilon }$
such that
\begin{equation}
\Psi \left( t_{\alpha },x_{\alpha },y_{\alpha }\right) >\sup_{\Sigma
_{\epsilon }}\Psi \left( t,x,y\right) -\alpha .  \label{lowb}
\end{equation}%
Set%
\begin{equation*}
\Psi _{\alpha }\left( t,x,y\right) =\Psi \left( t,x,y\right) -\frac{1}{%
2\alpha }\left[ \left\vert t-t_{\alpha }\right\vert ^{2}+\left\vert
x-x_{\alpha }\right\vert ^{2}+\left\vert y-y_{\alpha }\right\vert ^{2}\right]
.
\end{equation*}%
Clearly, $\Psi _{\alpha }\left( t,x,y\right) $ attains the maximum on $%
\Sigma _{\epsilon }$ at certain point $\left( t^{\ast },x^{\ast },y^{\ast
}\right) ,$ which depends implicitly on $\alpha $ and $\beta .$ One can
immediately get
\begin{equation}
\left\vert t^{\ast }-t_{\alpha }\right\vert ^{2}+\left\vert x^{\ast
}-x_{\alpha }\right\vert ^{2}+\left\vert y^{\ast }-y_{\alpha }\right\vert
^{2}\leq 2.  \label{bounded}
\end{equation}%
By the method applied in \cite{FS}, we can show that, for any given $%
\epsilon >0,$
\begin{equation}
\left\vert x^{\ast }-y^{\ast }\right\vert ^{2}=o\left( \alpha \right) .
\label{small}
\end{equation}%
Now we assume that the maximum of the function $\Psi _{\alpha }\left(
t,x,y\right) $ is always at the point $\left( T,x^{\ast },y^{\ast }\right) .$
Then, we have%
\begin{equation*}
\Psi _{\alpha }\left( T,x,y\right) \leq \Psi _{\alpha }\left( T,x^{\ast
},y^{\ast }\right) ,\text{ }\forall x,y\in \mathbb{R}^{n}.
\end{equation*}%
Thus
\begin{eqnarray*}
\underset{\alpha \rightarrow 0}{\text{limsup}}\left( u_{1}^{\epsilon }\left(
T,x^{\ast }\right) -u_{2}\left( T,y^{\ast }\right) \right) &\leq &\underset{%
\alpha \rightarrow 0}{\text{limsup}}\left( \left( 1-\epsilon \right)
u_{1}\left( T,x^{\ast }\right) -u_{2}\left( T,y^{\ast }\right) \right) \\
&\leq &\left( 1-\epsilon \right) \underset{\alpha \rightarrow 0}{\text{limsup%
}}\left( \Phi \left( x^{\ast }\right) -\Phi \left( y^{\ast }\right) \right)
\\
&&-\epsilon \underset{\alpha \rightarrow 0}{\text{limsup}}\Phi \left(
y^{\ast }\right) \\
&\leq &\epsilon C,
\end{eqnarray*}%
since $\Phi $ satisfies Lipschitz condition and (\ref{small}). The constant $%
C$ will change line to line. Next, we deal with%
\begin{eqnarray*}
&&u_{1}^{\epsilon }\left( t,x\right) -u_{2}\left( t,y\right) +\beta \left(
t-T\right) \\
&=&\Psi \left( t,x,x\right) \\
&\leq &\underset{\epsilon <t<T,x,y\in \mathbb{R}^{n}}{\text{limsup}}\Psi
\left( t,x,y\right) \\
&=&\lim_{\alpha \rightarrow 0}\Psi _{\alpha }\left( t_{\alpha },x_{\alpha
},y_{\alpha }\right) \\
&\leq &\underset{\alpha \rightarrow 0}{\text{limsup}}\Psi _{\alpha }\left(
T,x^{\ast },y^{\ast }\right) \\
&\leq &\underset{\alpha \rightarrow 0}{\text{limsup}}\left( u_{1}^{\epsilon
}\left( T,x^{\ast }\right) -u_{2}\left( T,y^{\ast }\right) \right) \\
&\leq &\epsilon C.
\end{eqnarray*}%
Here we have used the fact (\ref{lowb}). Hence, we have by the definition of
$u_{1}^{\epsilon }$%
\begin{eqnarray*}
u_{1}\left( t,x\right) -u_{2}\left( t,x\right) &=&\frac{1}{1-\epsilon }%
u_{1}^{\epsilon }\left( t,x\right) -u_{2}\left( t,x\right) +\frac{\epsilon }{%
\left( 1-\epsilon \right) \left( t-\epsilon \right) } \\
&=&u_{1}^{\epsilon }\left( t,x\right) -u_{2}\left( t,x\right) +\frac{%
\epsilon }{1-\epsilon }u_{1}^{\epsilon }\left( t,x\right) +\frac{\epsilon }{%
\left( 1-\epsilon \right) \left( t-\epsilon \right) } \\
&\leq &\epsilon C+\frac{\epsilon }{1-\epsilon }u_{1}^{\epsilon }\left(
t,x\right) +\frac{\epsilon }{\left( 1-\epsilon \right) \left( t-\epsilon
\right) }+\beta \left( T-t\right) .
\end{eqnarray*}%
Letting $\epsilon ,\beta \rightarrow 0,$ we obtain
\begin{equation}
u_{1}\left( t,x\right) \leq u_{2}\left( t,x\right) ,\text{ }\left(
t,x\right) \in \left[ 0,T\right] \times \mathbb{R}^{n}.  \label{r1}
\end{equation}%
From the arbitrary of $\left( t,x\right) ,$ we conclude that (\ref{r1})
holds on $\Sigma .$ Repeating the method above, we are able to get
\begin{equation}
u_{2}\left( t,x\right) \leq u_{1}\left( t,x\right) ,\text{ }\left(
t,x\right) \in \left[ 0,T\right] \times \mathbb{R}^{n}.  \label{r2}
\end{equation}%
Apparently, the rest proof will focus on the observation $t^{\ast }=T$
whenever $\epsilon ,\alpha $ are small enough and $\beta >0.$ We will prove this result by a contradiction. Suppose that $t^{\ast }<T.$

Define
\begin{equation*}
\phi \left( t,x,y\right) =\frac{1}{2\alpha }\left\vert x-y\right\vert ^{2}+%
\frac{\alpha }{2}\left[ \left\vert t-t_{\alpha }\right\vert ^{2}+\left\vert
x-x_{\alpha }\right\vert ^{2}+\left\vert y-y_{\alpha }\right\vert ^{2}\right]
-\beta T\left( t-T\right) .
\end{equation*}%
By virtue of Theorem 8.3 in \cite{CIL}, there exist $c_{i}\in \mathbb{R},$ $%
X_{i}\in \mathbb{S}^{n},$ $i=1,2$ such that%
\begin{eqnarray}
(\text{i})\text{ }\left( c_{1},D_{x_{1}}\phi \left( t^{\ast },x^{\ast
},y^{\ast }\right) ,X_{1}\right) &\in &\mathcal{P}^{2,+}u_{1}\left(
t,x\right) ,  \notag \\
\text{ }\left( c_{2},D_{x_{2}}\phi \left( t^{\ast },x^{\ast },y^{\ast
}\right) ,X_{2}\right) &\in &\mathcal{P}^{2,+}\left( -u_{2}\left( t,x\right)
\right),  \label{use1} \\
(\text{ii})\text{ }-\left( \frac{1}{\epsilon }+\left\Vert A\right\Vert
\right) I &\leq &\left(
\begin{array}{cc}
X_{1} & 0 \\
0 & X_{2}%
\end{array}%
\right) \leq A+\epsilon A^{2},  \label{use2} \\
(\text{iii})\text{ }c_{1}+c_{2} &=&\phi _{t}\left( t^{\ast },x^{\ast
},y^{\ast }\right) ,  \label{use3}
\end{eqnarray}%
where $A=D_{x}^{2}\phi \left( t^{\ast },x^{\ast },y^{\ast }\right) $ with
appropriate norm\footnote{%
Here, the norm of the symmetric matrix $A$ is
\begin{eqnarray*}
\left\Vert A\right\Vert &=&\sup \left\{ \left\vert \lambda \right\vert
:\lambda \text{ is an eigenvalue of }A\right\} \\
&=&\sup \left\{ \left\vert \left\langle A\xi ,\xi \right\rangle \right\vert
:\left\vert \xi \right\vert \leq 1\right\} ,
\end{eqnarray*}%
and $I$ denotes the identity matrix.}$.$ From the definition of $\phi $, it
is easy to show that%
\begin{eqnarray}
\phi _{t}\left( t^{\ast },x^{\ast },y^{\ast }\right) &=&-T\beta +\alpha
\left( t^{\ast }-t_{\alpha }\right) ,  \label{f1} \\
D_{x}\phi \left( t^{\ast },x^{\ast },y^{\ast }\right) &=&\frac{1}{\alpha }%
\left( x^{\ast }-y^{\ast }\right) +\alpha \left( x^{\ast }-x_{\alpha
}\right) ,  \label{f2} \\
D_{y}\phi \left( t^{\ast },x^{\ast },y^{\ast }\right) &=&-\frac{1}{\alpha }%
\left( x^{\ast }-y^{\ast }\right) +\alpha \left( y^{\ast }-y_{\alpha
}\right) ,  \label{f3} \\
D_{xx}^{2}\phi \left( t^{\ast },x^{\ast },y^{\ast }\right) &=&\left( \frac{1%
}{\alpha }+\alpha \right) I,\text{ }D_{xy}^{2}\phi \left( t^{\ast },x^{\ast
},y^{\ast }\right) =-\frac{1}{\alpha }I,  \label{f4} \\
D_{yy}^{2}\phi \left( t^{\ast },x^{\ast },y^{\ast }\right) &=&\left( \frac{1%
}{\alpha }+\alpha \right) I.  \label{f5}
\end{eqnarray}%
Clearly, (iii) in (\ref{use3}) can be expressed as
\begin{equation}
c_{1}+c_{2}=-T\beta +\alpha \left( t^{\ast }-t_{\alpha }\right)  \label{f6}
\end{equation}%
and
\begin{eqnarray}
A &=&\left(
\begin{array}{cc}
\left( \frac{1}{\alpha }+\alpha \right) I & -\frac{1}{\alpha }I \\
-\frac{1}{\alpha }I & \left( \frac{1}{\alpha }+\alpha \right) I%
\end{array}%
\right)  \notag \\
&=&\frac{1}{\alpha }\left(
\begin{array}{cc}
I & -I \\
-I & I%
\end{array}%
\right) +\alpha \left(
\begin{array}{cc}
I & 0 \\
0 & I%
\end{array}%
\right) .  \label{f7}
\end{eqnarray}%
Note that $\mathcal{P}^{2,+}\left( -u_{2}\left( t,x\right) \right) =-%
\mathcal{P}^{2,-}\left( u_{2}\left( t,x\right) \right) $ and for any $%
\lambda >0,$ we have%
\begin{equation*}
\mathcal{P}^{2,+}\left( \lambda u_{2}\right) =\lambda \mathcal{P}%
^{2,+}\left( u_{2}\right) ,\text{ }\mathcal{P}^{2,-}\left( \lambda
u_{2}\right) =\lambda \mathcal{P}^{2,-}\left( u_{2}\right) .
\end{equation*}%
Then, (i) in (\ref{use1}) can be rewritten as
\begin{eqnarray*}
\left( \frac{1}{1-\epsilon }c_{1},\frac{1}{1-\epsilon }D_{x}\phi \left(
t^{\ast },x^{\ast },y^{\ast }\right) ,\frac{1}{1-\epsilon }X_{1}\right) &\in
&\frac{1}{1-\epsilon }\mathcal{P}^{2,+}u_{1}^{\epsilon }\left( t,x\right) ,
\\
\left( -c_{2},-D_{y}\phi \left( t^{\ast },x^{\ast },y^{\ast }\right)
,-X_{2}\right) &\in &\mathcal{P}^{2,-}u_{2}\left( t,x\right) .
\end{eqnarray*}%
From the Definition \ref{dvis1}, it yields
\begin{equation}
\begin{array}{ll}
\min \Big (\frac{1}{1-\epsilon }D_{x}\phi \left( t^{\ast },x^{\ast },y^{\ast
}\right) G+K &  \\
\frac{1}{1-\epsilon }c_{1}+\mathcal{H}\left( t^{\ast },x^{\ast },\frac{1}{%
1-\epsilon }D_{x}\phi \left( t^{\ast },x^{\ast },y^{\ast }\right) ,\frac{1}{%
1-\epsilon }X_{1}\right) \Big ) & \geq 0%
\end{array}
\label{visform1}
\end{equation}%
and%
\begin{equation}
\begin{array}{ll}
\min \Big (-D_{y}\phi \left( t^{\ast },x^{\ast },y^{\ast }\right) G+K, &  \\
-c_{2}+\mathcal{H}\left( t^{\ast },x^{\ast },-D_{y}\phi \left( t^{\ast
},x^{\ast },y^{\ast }\right) ,-X_{2}\right) \Big ) & \leq 0.%
\end{array}
\label{visform2}
\end{equation}%
Obviously, from (\ref{visform1}) we always have
\begin{equation}
\left\{
\begin{array}{l}
\frac{1}{1-\epsilon }D_{x}\phi \left( t^{\ast },x^{\ast },y^{\ast }\right)
G+K\geq 0, \\
\frac{1}{1-\epsilon }c_{1}+\mathcal{H}\left( t^{\ast },x^{\ast },\frac{1}{%
1-\epsilon }D_{x}\phi \left( t^{\ast },x^{\ast },y^{\ast }\right) ,\frac{1}{%
1-\epsilon }X_{1}\right) \geq 0.%
\end{array}%
\right.  \label{visform11}
\end{equation}%
We will show that
\begin{equation}
-D_{y}\phi \left( t^{\ast },x^{\ast },y^{\ast }\right) G+K>0.  \label{uncon}
\end{equation}%
Assume that (\ref{uncon}) is not true, then for $1\leq i\leq m,$%
\begin{equation}
-\left( D_{y}\phi \left( t^{\ast },x^{\ast },y^{\ast }\right) G\right)
^{i}+K^{i}<0,  \label{uncon2}
\end{equation}%
where $\left( \eta \right) ^{i}$ stands for the $i$th coordinate of $\eta $
when $\eta $ is row vector. Now multiply the $i$th inequality of the first
term in (\ref{visform1}) by $1-\epsilon $ to get
\begin{equation}
\left( D_{x}\phi \left( t^{\ast },x^{\ast },y^{\ast }\right) G\right)
^{i}+\left( 1-\epsilon \right) K^{i}>0.  \label{uncon3}
\end{equation}%
From (\ref{uncon2}), (\ref{uncon3}), (\ref{f2}) and (\ref{f3}), we have%
\begin{eqnarray*}
\epsilon K^{i} &<&\left( D_{x}\phi \left( t^{\ast },x^{\ast },y^{\ast
}\right) ^{\top }G\right) ^{i}+\left( D_{y}\phi \left( t^{\ast },x^{\ast
},y^{\ast }\right) ^{\top }G\right) ^{i} \\
&=&\frac{1}{\alpha }\left( x^{\ast }-y^{\ast }\right) +\alpha \left( x^{\ast
}-x_{\alpha }\right) -\frac{1}{\alpha }\left( x^{\ast }-y^{\ast }\right)
+\alpha \left( y^{\ast }-y_{\alpha }\right) \\
&=&\alpha \left( \left[ \left( x^{\ast }-x_{\alpha }\right) ^{\top }+\left(
y^{\ast }-y_{\alpha }\right) ^{\top }\right] G\right) ^{i}.
\end{eqnarray*}%
So when we fix $\epsilon ,$ and take $\alpha $ small enough, it will lead a
contradiction since $K^{i}>k_{0}>0,$ $1\leq i\leq m.$ Therefore, when $%
\alpha $ is small, (\ref{visform2}) is equivalent to
\begin{equation}
-c_{2}+\mathcal{H}\left( t^{\ast },x^{\ast },-D_{y}\phi \left( t^{\ast
},x^{\ast },y^{\ast }\right) ,-X_{2}\right) >0.  \label{visform3}
\end{equation}%
Combing (\ref{visform11}) and (\ref{visform3}), we get%
\begin{eqnarray}
&&c_{1}+\left( 1-\epsilon \right) \mathcal{H}\left( t^{\ast },x^{\ast },%
\frac{1}{1-\epsilon }D_{x}\phi \left( t^{\ast },x^{\ast },y^{\ast }\right) ,%
\frac{1}{1-\epsilon }X_{1}\right)  \notag \\
&&+c_{2}-\mathcal{H}\left( t^{\ast },x^{\ast },-D_{y}\phi \left( t^{\ast
},x^{\ast },y^{\ast }\right) ,-X_{2}\right)  \notag \\
&>&0.  \label{contrdi1}
\end{eqnarray}%
Noting (\ref{f1}), (\ref{f6}) and (\ref{contrdi1}) can be expressed as
\begin{eqnarray}
&&\left( 1-\epsilon \right) \mathcal{H}\left( t^{\ast },x^{\ast },\frac{1}{%
1-\epsilon }D_{x}\phi \left( t^{\ast },x^{\ast },y^{\ast }\right) ,\frac{1}{%
1-\epsilon }X_{1}\right)  \notag \\
&&-\mathcal{H}\left( t^{\ast },x^{\ast },-D_{y}\phi \left( t^{\ast },x^{\ast
},y^{\ast }\right) ,-X_{2}\right)  \notag \\
&>&\beta T-\alpha \left( t^{\ast }-t_{\alpha }\right) .  \label{contrdi2}
\end{eqnarray}%
Then, we will estimate (\ref{contrdi2}) step by step. At the beginning,
\begin{eqnarray*}
&&\left( 1-\epsilon \right) \mathcal{H}\left( t^{\ast },x^{\ast },\frac{1}{%
1-\epsilon }D_{x}\phi \left( t^{\ast },x^{\ast },y^{\ast }\right) ,\frac{1}{%
1-\epsilon }X_{1}\right) \\
&=&\left( 1-\epsilon \right) \inf_{v\in U}\Bigg \{\frac{1}{2}\text{\textrm{Tr%
}}\left( \sigma \sigma ^{\ast }\left( t^{\ast },x^{\ast },v\right) \frac{1}{%
1-\epsilon }X_{1}\right) \\
&&+\frac{1}{1-\epsilon }\left\langle D_{x}\phi \left( t^{\ast },x^{\ast
},y^{\ast }\right) ,b\left( t^{\ast },x^{\ast },v\right) \right\rangle \\
&&+f\left( t^{\ast },x^{\ast },u_{1}^{\epsilon }\left( t^{\ast },x^{\ast
}\right) ,\frac{1}{1-\epsilon }D_{x}\phi \left( t^{\ast },x^{\ast },y^{\ast
}\right) \sigma \left( t^{\ast },x^{\ast },v\right) ,v\right) \Bigg \},
\end{eqnarray*}%
whilst%
\begin{eqnarray*}
&&\mathcal{H}\left( t^{\ast },y^{\ast },-D_{y}\phi \left( t^{\ast },x^{\ast
},y^{\ast }\right) ,-X_{2}\right) \\
&=&\inf_{v\in U}\Bigg \{-\frac{1}{2}\text{\textrm{Tr}}\left( \sigma \sigma
^{\ast }\left( t^{\ast },y^{\ast },v\right) X_{2}\right) +\left\langle
-D_{y}\phi \left( t^{\ast },x^{\ast },y^{\ast }\right) ,b\left( t^{\ast
},y^{\ast },v\right) \right\rangle \\
&&+f\left( t^{\ast },y^{\ast },u_{2}\left( t^{\ast },y^{\ast }\right)
,-D_{y}\phi \left( t^{\ast },x^{\ast },y^{\ast }\right) \sigma \left(
t^{\ast },y^{\ast },v\right) ,v\right) \Bigg \}.
\end{eqnarray*}%
Next, (\ref{contrdi2}) can be rewritten as
\begin{eqnarray*}
&&\left( 1-\epsilon \right) \mathcal{H}\left( t^{\ast },x^{\ast },\frac{1}{%
1-\epsilon }D_{x}\phi \left( t^{\ast },x^{\ast },y^{\ast }\right) ,\frac{1}{%
1-\epsilon }X_{1}\right) \\
&&-\mathcal{H}\left( t^{\ast },x^{\ast },-D_{y}\phi \left( t^{\ast },x^{\ast
},y^{\ast }\right) ,-X_{2}\right) \\
&\leq &\sup_{v\in U}\Bigg \{\frac{1}{2}\text{\textrm{Tr}}\left( \sigma
\sigma ^{\ast }\left( t^{\ast },x^{\ast },v\right) X_{1}\right) +\frac{1}{2}%
\text{\textrm{Tr}}\left( \sigma \sigma ^{\ast }\left( t^{\ast },y^{\ast
},v\right) X_{2}\right) \\
&&+\left\langle D_{x}\phi \left( t^{\ast },x^{\ast },y^{\ast }\right)
,b\left( t^{\ast },x^{\ast },v\right) \right\rangle +\left\langle D_{y}\phi
\left( t^{\ast },x^{\ast },y^{\ast }\right) ,b\left( t^{\ast },y^{\ast
},v\right) \right\rangle \\
&&+\left( 1-\epsilon \right) f\left( t^{\ast },x^{\ast },u_{1}^{\epsilon
}\left( t^{\ast },x^{\ast }\right) ,\frac{1}{1-\epsilon }D_{x}\phi \left(
t^{\ast },x^{\ast },y^{\ast }\right) \sigma \left( t^{\ast },x^{\ast
},v\right) ,v\right) \\
&&-f\left( t^{\ast },y^{\ast },u_{2}\left( t^{\ast },y^{\ast }\right)
,-D_{y}\phi \left( t^{\ast },x^{\ast },y^{\ast }\right) \sigma \left(
t^{\ast },y^{\ast },v\right) ,v\right) \Bigg \}.
\end{eqnarray*}%
Set
\begin{eqnarray*}
\Upsilon _{1} &=&\frac{1}{2}\text{\textrm{Tr}}\left( \sigma \sigma ^{\ast
}\left( t^{\ast },x^{\ast },v\right) X_{1}\right) +\frac{1}{2}\text{\textrm{%
Tr}}\left( \sigma \sigma ^{\ast }\left( t^{\ast },y^{\ast },v\right)
X_{2}\right) , \\
\Upsilon _{2} &=&\left\langle D_{x}\phi \left( t^{\ast },x^{\ast },y^{\ast
}\right) ,b\left( t^{\ast },x^{\ast },v\right) \right\rangle +\left\langle
D_{y}\phi \left( t^{\ast },x^{\ast },y^{\ast }\right) ,b\left( t^{\ast
},y^{\ast },v\right) \right\rangle , \\
\Upsilon _{3} &=&\left( 1-\epsilon \right) f\left( t^{\ast },x^{\ast
},u_{1}^{\epsilon }\left( t^{\ast },x^{\ast }\right) ,\frac{1}{1-\epsilon }%
D_{x}\phi \left( t^{\ast },x^{\ast },y^{\ast }\right) \sigma \left( t^{\ast
},x^{\ast },v\right) ,v\right) \\
&&-f\left( t^{\ast },y^{\ast },u_{2}\left( t^{\ast },y^{\ast }\right)
,-D_{y}\phi \left( t^{\ast },x^{\ast },y^{\ast }\right) \sigma \left(
t^{\ast },y^{\ast },v\right) ,v\right) .
\end{eqnarray*}%
For $\Upsilon _{1},$ with $\sigma _{1}=\sigma \left( t^{\ast },x^{\ast
},v\right) $, $\sigma _{2}=\sigma \left( t^{\ast },y^{\ast },v\right) ,$
from (\ref{use2}) and (\ref{f7}), we have%
\begin{eqnarray}
&&\frac{1}{2}\text{\textrm{Tr}}\left( \sigma \sigma ^{\ast }\left( t^{\ast
},x^{\ast },v\right) X_{1}\right) +\frac{1}{2}\text{\textrm{Tr}}\left(
\sigma \sigma ^{\ast }\left( t^{\ast },y^{\ast },v\right) X_{2}\right)
\notag \\
&=&\frac{1}{2}\text{\textrm{Tr}}\left(
\begin{array}{cc}
\sigma _{1}\sigma _{1}^{\top } & \sigma _{1}\sigma _{2}^{\top } \\
\sigma _{2}\sigma _{1}^{\top } & \sigma _{2}\sigma _{2}^{\top }%
\end{array}%
\right) \left(
\begin{array}{cc}
X_{1} & 0 \\
0 & X_{2}%
\end{array}%
\right)  \notag \\
&\leq &\frac{1}{2}\left( \frac{3}{\alpha }+2\alpha \right) \left\Vert \sigma
_{1}-\sigma _{1}\right\Vert ^{2}+C\frac{1}{2}\left( \alpha +\alpha
^{3}\right)  \notag \\
&\leq &\frac{C^{2}}{2}\left( \frac{3}{\alpha }+2\alpha \right) \left\vert
x^{\ast }-y^{\ast }\right\vert ^{2}+C\frac{1}{2}\left( \alpha +\alpha
^{3}\right)  \notag \\
&\leq &o\left( 1\right) ,  \label{fr1}
\end{eqnarray}%
where we have used the the assumption that $\sigma $ is bounded and
Lipschitz condition.

For $\Upsilon _{2},$ we derive%
\begin{eqnarray}
&&\left\langle D_{x}\phi \left( t^{\ast },x^{\ast },y^{\ast }\right)
,b\left( t^{\ast },x^{\ast },v\right) \right\rangle +\left\langle D_{y}\phi
\left( t^{\ast },x^{\ast },y^{\ast }\right) ,b\left( t^{\ast },y^{\ast
},v\right) \right\rangle  \notag \\
&=&\left\langle \frac{1}{\alpha }\left( x^{\ast }-y^{\ast }\right) +\alpha
\left( x^{\ast }-x_{\alpha }\right) ,b\left( t^{\ast },x^{\ast },v\right)
\right\rangle  \notag \\
&&+\left\langle -\frac{1}{\alpha }\left( x^{\ast }-y^{\ast }\right) +\alpha
\left( y^{\ast }-y_{\alpha }\right) ,b\left( t^{\ast },y^{\ast },v\right)
\right\rangle  \notag \\
&\leq &\frac{C}{\alpha }\left\vert x^{\ast }-y^{\ast }\right\vert
^{2}+2C\alpha  \notag \\
&\leq &o\left( 1\right) .  \label{fr2}
\end{eqnarray}

For $\Upsilon _{3},$ we handle%
\begin{eqnarray*}
&&\left( 1-\epsilon \right) f\left( t^{\ast },x^{\ast },u_{1}^{\epsilon
}\left( t^{\ast },x^{\ast }\right) ,\frac{1}{1-\epsilon }D_{x}\phi \left(
t^{\ast },x^{\ast },y^{\ast }\right) \sigma \left( t^{\ast },x^{\ast
},v\right) ,v\right) \\
&&-f\left( t^{\ast },y^{\ast },u_{2}\left( t^{\ast },y^{\ast }\right)
,-D_{y}\phi \left( t^{\ast },x^{\ast },y^{\ast }\right) \sigma \left(
t^{\ast },y^{\ast },v\right) \right) \\
&\leq &\left( 1-\epsilon \right) \Bigg [\varpi \left( \left\vert x^{\ast
}-y^{\ast }\right\vert \left( 1+\frac{1}{1-\epsilon }\right) D_{x}\phi
\left( t^{\ast },x^{\ast },y^{\ast }\right) \sigma \left( t^{\ast },x^{\ast
},v\right) \right) \\
&&+C\left\vert u_{1}^{\epsilon }\left( t^{\ast },x^{\ast }\right)
-u_{2}\left( t^{\ast },y^{\ast }\right) \right\vert \\
&&+C\left\vert \frac{1}{1-\epsilon }D_{x}\phi \left( t^{\ast },x^{\ast
},y^{\ast }\right) \sigma \left( t^{\ast },x^{\ast },v\right) +D_{y}\phi
\left( t^{\ast },x^{\ast },y^{\ast }\right) \sigma \left( t^{\ast },y^{\ast
},v\right) \right\vert ,
\end{eqnarray*}%
where we have used the Lipschitz continuity of $f.$ Now%
\begin{eqnarray*}
&&\left\vert \frac{1}{1-\epsilon }D_{x}\phi \left( t^{\ast },x^{\ast
},y^{\ast }\right) \sigma \left( t^{\ast },x^{\ast },v\right) +D_{y}\phi
\left( t^{\ast },x^{\ast },y^{\ast }\right) \sigma \left( t^{\ast },y^{\ast
},v\right) \right\vert \\
&\leq &\frac{1}{\alpha \left( 1-\epsilon \right) }\left\vert \left( x^{\ast
}-y^{\ast }\right) \left( \sigma \left( t^{\ast },x^{\ast },v\right) -\sigma
\left( t^{\ast },y^{\ast },v\right) \right) \right\vert \\
&&+\frac{1}{1-\epsilon }\left\vert \alpha \left( x^{\ast }-x_{\alpha
}\right) \sigma \left( t^{\ast },x^{\ast },v\right) +\alpha \left( y^{\ast
}-y_{\alpha }\right) \sigma \left( t^{\ast },y^{\ast },v\right) \right\vert
\\
&\leq &\frac{\left\vert x^{\ast }-y^{\ast }\right\vert ^{2}}{\alpha \left(
1-\epsilon \right) }+\frac{2\alpha C}{1-\epsilon },
\end{eqnarray*}%
where we have used the bound of $\sigma $ and the fact (\ref{bounded}). So
\begin{eqnarray}
&&\left( 1-\epsilon \right) f\left( t^{\ast },x^{\ast },u_{1}^{\epsilon
}\left( t^{\ast },x^{\ast }\right) ,\frac{1}{1-\epsilon }D_{x}\phi \left(
t^{\ast },x^{\ast },y^{\ast }\right) \sigma \left( t^{\ast },x^{\ast
},v\right) ,v\right)  \notag \\
&&-f\left( t^{\ast },y^{\ast },u_{2}\left( t^{\ast },y^{\ast }\right)
,-D_{y}\phi \left( t^{\ast },x^{\ast },y^{\ast }\right) \sigma \left(
t^{\ast },y^{\ast },v\right) ,v\right)  \notag \\
&\leq &\left( 1-\epsilon \right) \varpi \left( \left\vert x^{\ast }-y^{\ast
}\right\vert \left( 1+\frac{1}{1-\epsilon }D_{x}\phi \left( t^{\ast
},x^{\ast },y^{\ast }\right) \sigma \left( t^{\ast },x^{\ast },v\right)
\right) \right)  \notag \\
&&+C\epsilon +\beta \left( 1-t^{\ast }\right) +\frac{\left\vert x^{\ast
}-y^{\ast }\right\vert ^{2}}{\alpha \left( 1-\epsilon \right) }+\frac{%
2\alpha C}{1-\epsilon }.  \label{fr3}
\end{eqnarray}%
Now from (\ref{fr1}), (\ref{fr2}) and (\ref{fr3}), we have
\begin{eqnarray*}
&&\beta T-\alpha \left( t^{\ast }-t_{\alpha }\right) \\
&\leq &o\left( 1\right) +\left( 1-\epsilon \right) \varpi \left( \left\vert
x^{\ast }-y^{\ast }\right\vert \left( 1+\frac{1}{1-\epsilon }D_{x}\phi
\left( t^{\ast },x^{\ast },y^{\ast }\right) \sigma \left( t^{\ast },x^{\ast
},v\right) \right) \right) \\
&&+C\epsilon +\beta \left( T-t^{\ast }\right) +\frac{\left\vert x^{\ast
}-y^{\ast }\right\vert ^{2}}{\alpha \left( 1-\epsilon \right) }+\frac{%
2\alpha C}{1-\epsilon }.
\end{eqnarray*}%
Letting $\epsilon ,$ $\alpha \rightarrow 0,$ and noticing $t^{\ast }<T,$ we
obtain%
\begin{equation*}
\beta t^{\ast }\leq 0,
\end{equation*}%
which leads to a contradiction. Consequently, $t^{\ast }=T$ holds true. We
thus complete the proof.\hfill $\Box $

\begin{remark}
\label{existence}The present paper considers a stochastic optimal singular
control problem, in which the cost function is defined through a backward
stochastic differential equation with singular control. It has been shown
that the value function is a unique viscosity solution of the corresponding
Hamilton-Jacobi-Bellman inequality (\ref{pde1}), in a given class of
continuous and bounded functions. Nonetheless, as an important issue, the
\emph{existence} of optimal singular controls has never been exploited. In
Haussmannand, Suo \cite{HS1}, the authors apply the compactification method
to study the classical and singular control problem of It\^{o}'s type of
stochastic differential equation, where the problem is reformulated as a
martingale problem on an appropriate canonical space after the relaxed form
of the classical control is introduced. Under some mild continuity
assumptions on the data, they obtain the existence of optimal control by
purely probabilistic arguments. Nevertheless, in the framework of BSDE, the
trajectory of $Y$ seems to be a c\`{a}dl\`{a}g process (from French, for
right continuous with left hand limits). Hence, we may consider $Y$ in some
space with appropriate topologies, for instance, Meyer-Zheng topology and
obtain the convergence of probability measures deduced by $Y$ involving
relaxed control. For this topic, reader can also refer articles \cite{HS1,
FH} in this direction. We point out that the related work from the technique
of PDEs (see \cite{BLRT, BGM} therein). From Wang \cite{W}, it is possible
to construct the optimal control via the existence of diffusion with
refections (see \cite{CMR}). However, it is interesting to expand this
result to forward backward stochastic systems. We shall consider these
crucial topics in future. As this complete remake of the existence is much
longer than the present paper, it will be reported elsewhere.
\end{remark}


\section{Example}

\label{sect:4}

In this section, we illustrate our theoretic result by looking at a concrete
example.

\begin{example}
Consider the following controlled SDE $\left( n=1\right) ,$%
\begin{equation*}
\left\{
\begin{array}{lll}
\mathrm{d}X_{s}^{t,x;v,\xi } & = & \left[ X_{s}^{t,x;v,\xi
}+X_{s}^{t,x;v,\xi }v_{s}\right] \mathrm{d}s+X_{s}^{t,x;v,\xi }v_{s}\mathrm{d%
}W_{s} \\
&  & +G\mathrm{d}\xi _{s},\text{ }s\in \left[ t,T\right] , \\
X_{t}^{t,x;v,\xi } & = & x,%
\end{array}%
\right.
\end{equation*}%
with the control domain $U=\left[ -1,0\right] \cup \left[ 1,2\right] $,
constants $G>0$ and $K>0$.

The cost functional is defined by%
\begin{equation*}
\left\{
\begin{array}{lll}
\mathrm{d}Y_{s}^{t,x;v,\xi } & = & Z_{s}^{t,x;v,\xi }v_{s}\mathrm{d}%
s+Z_{s}^{t,x;v,\xi }\mathrm{d}W_{s}+K\mathrm{d}\xi _{s},\text{ }s\in \left[
t,T\right] , \\
Y_{T}^{t,x;v,\xi } & = & X_{T}^{t,x;v,\xi }.%
\end{array}%
\right.
\end{equation*}%
It is easy to give the associated H-J-B inequality is
\begin{equation}
\left\{
\begin{array}{l}
\min \left\{ G_{t}u_{x}+K,u_{t}+xu_{x}+\min_{v\in U}\left[ \frac{1}{2}%
u_{xx}x^{2}v^{2}+xvu_{x}-xv^{2}u_{x}\right] \right\} =0, \\
u\left( T,x\right) =x,\text{ for all }x\in \mathbb{R}^{n}.%
\end{array}%
\right.  \label{exa}
\end{equation}%
By Theorem \ref{t4}, we deal with the second term in (\ref{exa}), then
obtain
\begin{equation*}
u\left( t,x\right) =\left\{
\begin{array}{c}
e^{t-T}x,\text{ if }x>0, \\
e^{T-t}x,\text{ if }x\leq 0.%
\end{array}%
\right.
\end{equation*}%
Assume that the optimal control is $\left( v^{\ast },\xi ^{\ast }\right) $.
If $x\leq 0$, then $u\left( s,X^{t,x;v^{\ast },\xi ^{\ast }}\left( s\right)
\right) =e^{T-s}X^{t,x;v^{\ast },\xi ^{\ast }}\left( s\right) $, which
yields that%
\begin{equation*}
\mathcal{P}^{2,+}u\left( s,X^{t,x;v^{\ast },\xi ^{\ast }}\left( s\right)
\right) =\left[ e^{t-T}x,+\infty \right) \times \left[ e^{t-T},+\infty
\right) \times \left[ 0,+\infty \right) ,
\end{equation*}%
and
\begin{equation*}
\mathcal{P}^{2,-}u\left( s,X^{t,x;v^{\ast },\xi ^{\ast }}\left( s\right)
\right) =\left( -\infty ,e^{t-T}x\right] \times \left( -\infty ,e^{t-T}%
\right] \times \left( -\infty ,0\right] .
\end{equation*}%
It is fairly easy to verify that, for any $\left( p,q,X\right) \in \mathcal{P%
}^{2,+}u\left( t,x\right) $
\begin{eqnarray*}
Gq+K &\geq &0, \\
p+\frac{1}{2}Xx^{2}-xq &\geq &0,
\end{eqnarray*}%
which means that $u\left( t,x\right) $ is the subsolution when $x\leq 0$.
Moreover, for $\left( \tilde{p},\tilde{q},\tilde{X}\right) \in \mathcal{P}%
^{2,-}u\left( t,x\right) ,$ we have%
\begin{eqnarray*}
G\tilde{q}+K &\geq &0, \\
\tilde{p}+\frac{1}{2}\tilde{X}x^{2}-x\tilde{q} &\leq &0,
\end{eqnarray*}%
which means $u\left( t,x\right) $ is the supersolution.
\end{example}

\appendix\label{sect:5}

\section{Forward-backward SDEs with Singular Controls}

In this section, we revisit the theory of FBSDEs with controls. We will give
some useful estimates for (\ref{eq1}). Consider the following SDE with
controls by the initial condition $\left( t,x\right) \in \left[ 0,T\right]
\times \mathbb{R}^{n}$,%
\begin{equation}
\left\{
\begin{array}{rcl}
\mathrm{d}X_{s}^{t,x,v,\xi } & = & b\left( s,X_{s}^{t,x,v,\xi },v_{s}\right)
\mathrm{d}s+\sigma \left( s,X_{s}^{t,x,v,\xi },v_{s}\right) \mathrm{d}%
W_{s}+G_{s}\mathrm{d}\zeta _{s}, \\
X_{t}^{t,x,v,\xi } & = & x,\text{ }s\in \left[ 0,T\right] .%
\end{array}%
\right.  \label{ASDE}
\end{equation}

\begin{lemma}
\label{state1}Assume that \emph{(A1)-(A3)} hold, SDE \emph{(\ref{ASDE})}
admits a strong solution, and for any $t\in \left[ 0,T\right] ,$ and $x,$ $%
x^{\prime }\in \mathbb{R}^{n},$ there exists positive constant $C>0$ ($C$
depends only on the Lipschitz and linear growth constant for $b$, $\sigma ,$
with respect to $x$ and $u$, $G$ and $\xi $)$,$%
\begin{eqnarray}
\mathbb{E}\left[ \sup_{t\leq s\leq T}\left\vert X_{s}^{t,x,v,\xi
}-X_{s}^{t,x^{\prime },v,\xi }\right\vert ^{2}\right] &\leq &C\left\vert
x-x^{\prime }\right\vert ^{2},  \label{E11} \\
\mathbb{E}\left[ \sup_{t\leq s\leq T}\left\vert X_{s}^{t,x,v,\xi
}\right\vert ^{2}\right] &\leq &C\left( 1+\left\vert x\right\vert
^{2}\right) .  \label{E22}
\end{eqnarray}
\end{lemma}

\paragraph{Proof}

The estimate (\ref{E11}) can be seen in Theorem 4.1 in \cite{HS2}.
Inequality (\ref{E22}) can be obtained by the classical approach
(Burkholder--Davis--Gundy inequality and Gronwall's Lemma), so we omit the
proof. 
~\hfill $\Box $


\begin{lemma}
\label{state2}Assume that \emph{(A1)-(A3)} hold. Then for all $0\leq t\leq
T, $ $x,$ $x^{\prime }\in \mathbb{R}^{n},$ we have%
\begin{equation}
\mathbb{E}^{\mathscr{F}_{t}}\left[ \sup_{t\leq s\leq T}\left\vert
Y_{s}^{t,x,v,\xi }-Y_{s}^{t,x^{\prime },v,\xi }\right\vert
^{2}+\int_{s}^{T}\left\vert Z_{s}^{t,x,v,\xi }-Z_{s}^{t,x^{\prime },v,\xi
}\right\vert ^{2}\mathrm{d}s\right] \leq C\left\vert x-x^{\prime
}\right\vert ,\text{ }a.s.  \label{E1}
\end{equation}%
\begin{equation}
\mathbb{E}^{\mathscr{F}_{t}}\left[ \sup_{t\leq s\leq T}\left\vert
Y_{s}^{t,x,v,\xi }\right\vert ^{2}+\int_{s}^{T}\left\vert Z_{s}^{t,x,v,\xi
}\right\vert ^{2}\mathrm{d}s\right] \leq C\left( 1+\left\vert x\right\vert
^{2}\right) ,\text{ }a.s.  \label{E2}
\end{equation}%
\begin{equation}
\mathbb{E}^{\mathscr{F}_{t}}\left[ \left\vert X_{t}^{t,x,v,\xi
}-X_{t}^{t,x,v^{\prime },\xi ^{\prime }}\right\vert ^{2}\right] \leq C%
\mathbb{E}^{\mathscr{F}_{t}}\left[ \int_{s}^{T}\left\vert
v_{s}-v_{s}^{\prime }\right\vert ^{2}\mathrm{d}s+\left\vert \xi _{s}-\xi
_{s}^{\prime }\right\vert ^{2}\right] ,\text{ }a.s.  \label{EE}
\end{equation}%
\begin{eqnarray}
&&\mathbb{E}^{\mathscr{F}_{t}}\left[ \left\vert Y_{t}^{t,x,v,\xi
}-Y_{t}^{t,x,v^{\prime },\xi ^{\prime }}\right\vert
^{2}+\int_{t}^{T}\left\vert Z_{s}^{t,x,v,\xi }-Z_{s}^{t,x,v^{\prime },\xi
^{\prime }}\right\vert ^{2}\mathrm{d}s\right]  \notag \\
&\leq &C\mathbb{E}^{\mathscr{F}_{t}}\left[ \int_{s}^{T}\left( \left\vert
v_{s}-v_{s}^{\prime }\right\vert ^{2}+\left\vert \xi _{s}-\xi _{s}^{\prime
}\right\vert ^{2}\right) \mathrm{d}s+\left\vert \xi _{T}-\xi _{T}^{\prime
}\right\vert ^{2}\right] ,\text{ }a.s.  \label{EE4}
\end{eqnarray}
\end{lemma}

\paragraph{Proof}

Suppose $\left( Y_{s}^{t,x;v,\xi },Z_{s}^{t,x;v,\xi }\right) $ and $\left(
Y_{s}^{t,x^{\prime },v,\xi },Z_{s}^{t,x^{\prime },v,\xi }\right) $ are the
solutions to the following BSDEs, respectively,%
\begin{eqnarray}
Y_{s}^{t,x;v,\xi } &=&\Phi \left( X_{T}^{t,x;v,\xi }\right)
+\int_{s}^{T}f\left( r,X_{r}^{t,x;v,\xi },Y_{r}^{t,x;v,\xi
},Z_{r}^{t,x;v,\xi },v_{r}\right) \mathrm{d}s  \notag \\
&&-\int_{s}^{T}Z_{r}^{t,x;v,\xi }\mathrm{d}W_{r}+K\xi _{T}-K\xi _{s},
\label{B1}
\end{eqnarray}%
\begin{eqnarray}
Y_{s}^{t,x^{\prime };v,\xi } &=&\Phi \left( X_{T}^{t,x^{\prime };v,\xi
}\right) +\int_{s}^{T}f\left( r,X_{r}^{t,x^{\prime };v,\xi
},Y_{r}^{t,x^{\prime };v,\xi },Z_{r}^{t,x^{\prime };v,\xi },v_{r}\right)
\mathrm{d}s  \notag \\
&&-\int_{s}^{T}Z_{r}^{t,x^{\prime };v,\xi }\mathrm{d}W_{r}+K\xi _{T}-K\xi
_{s}.  \label{B2}
\end{eqnarray}%
We can get estimate (\ref{E1}) by B-D-G inequality and (\ref{E11}). Put $%
\bar{Y}_{s}^{t,x,v,\xi }=Y_{s}^{t,x;v,\xi }+K\xi _{s}.$ Then Eq. (\ref{B1})
becomes to
\begin{eqnarray*}
\bar{Y}_{s}^{t,x;v,\xi } &=&\bar{Y}_{T}^{t,x;v,\xi }+\int_{s}^{T}f\left(
t,X_{r}^{t,x;v,\xi },\bar{Y}_{r}^{t,x;v,\xi }-K\xi _{r},Z_{r}^{t,x;v,\xi
},v_{r}\right) \mathrm{d}r \\
&&-\int_{s}^{T}Z_{r}^{t,x;v,\xi }\mathrm{d}W_{r}.
\end{eqnarray*}%
Applying It\^{o}'s formula to $\left\vert \bar{Y}_{s}^{t,x,v,\xi
}\right\vert ^{2}$ on the internal $\left[ s,T\right] ,$ we have%
\begin{eqnarray*}
&&\mathbb{E}\left[ \left\vert Y_{s}^{t,x,v,\xi }\right\vert ^{2}+\frac{1}{2}%
\int_{s}^{T}\left\vert Z_{r}^{t,x;v,\xi }\right\vert ^{2}\mathrm{d}r\right]
\leq \mathbb{E}\left[ \left\vert \bar{Y}_{s}^{t,x,v,\xi }\right\vert ^{2}+%
\frac{1}{2}\int_{s}^{T}\left\vert Z_{r}^{t,x;v,\xi }\right\vert ^{2}\mathrm{d%
}r\right] \\
&=&\mathbb{E}^{\mathscr{F}_{t}}\left[ \left\vert \bar{Y}_{T}^{t,x,v,\xi
}\right\vert ^{2}+2\int_{s}^{T}\bar{Y}_{r}^{t,x;v,\xi }f\left(
t,X_{r}^{t,x;v,\xi },\bar{Y}_{r}^{t,x;v,\xi }-K\xi _{r},Z_{r}^{t,x;v,\xi
},v_{r}\right) \mathrm{d}r\right] \\
&\leq &C\mathbb{E}^{\mathscr{F}_{t}}\Bigg [1+\left\vert x\right\vert
^{2}+\left\vert K\xi _{T}\right\vert ^{2}+\int_{s}^{T}\left\vert
v_{r}\right\vert ^{2}\mathrm{d}r+\int_{s}^{T}\left\vert Y_{r}^{t,x;v,\xi
}\right\vert ^{2}\mathrm{d}r\Bigg ].
\end{eqnarray*}%
By the Grownwall inequality, we can get the estimate (\ref{E2}).

Set%
\begin{eqnarray*}
\tilde{Y}_{s}^{t,v,\xi } &:&=Y_{s}^{t,x,v,\xi }-Y_{s}^{t,x^{\prime },v,\xi },
\\
\tilde{Z}_{s}^{t,v,\xi } &:&=Z_{s}^{t,x,v,\xi }-Z_{s}^{t,x^{\prime },v,\xi },
\\
\tilde{f}_{r} &:&=f\left( r,X_{r}^{t,x;v,\xi },Y_{r}^{t,x;v,\xi
},Z_{r}^{t,x;v,\xi },v_{r}\right) \\
&&-f\left( r,X_{r}^{t,x^{\prime };v,\xi },Y_{r}^{t,x^{\prime };v,\xi
},Z_{r}^{t,x^{\prime };v,\xi },v_{r}\right) .
\end{eqnarray*}%
Then, consider the following BSDE:%
\begin{equation*}
\tilde{Y}_{s}^{t,v,\xi }=\tilde{Y}_{T}^{t,v,\xi }+\int_{s}^{T}\tilde{f}_{r}%
\mathrm{d}r-\int_{s}^{T}\tilde{Z}_{r}^{t,v,\xi }\mathrm{d}W_{r}.
\end{equation*}%
Using Ito's formula to $\left\vert X_{t}^{t,x,v,\xi }-X_{t}^{t,x,v^{\prime
},\xi ^{\prime }}\right\vert ^{2}$ and $\left\vert Y_{t}^{t,x,v,\xi
}-Y_{t}^{t,x,v^{\prime },\xi ^{\prime }}\right\vert ^{2}$, respectively, it
is easy to get (\ref{EE}) and (\ref{EE4}). We omit the proof of since its
similar to (\ref{E2}). ~\hfill $\Box $

\begin{lemma}
\label{state4}Consider the BSDEs \emph{(\ref{BV1})} and \emph{(\ref{BV2})}.
There exists a positive constant $C$ such that, for any $\left( v,\xi
\right) \in \mathcal{U}$, we have
\begin{equation}
\left\vert Y_{t}^{1,v,0}-Y_{t}^{2,v}\right\vert \leq C\delta ^{\frac{3}{2}},%
\text{ a.s.}
\end{equation}
\end{lemma}

\paragraph{Proof}

By classical estimate to SDE (\ref{ASDE}), we have for some positive
constant $C$%
\begin{equation*}
\mathbb{E}^{\mathscr{F}_{t}}\left[ \sup_{t\leq s\leq t+\delta }\left\vert
X^{t,x;v,0}-x\right\vert ^{2}\right] \leq C\delta \left( 1+\left\vert
x\right\vert ^{2}\right) .
\end{equation*}%
Applying the It\^{o}'s formula to BSDE (\ref{BV1}) and (\ref{BV2}), we have
\begin{eqnarray*}
&&\mathbb{E}^{\mathscr{F}_{t}}\left[ \int_{t}^{t+\delta }\left( \left\vert
Y_{s}^{1,v,0}-Y_{s}^{2,v}\right\vert ^{2}+\left\vert
Z_{s}^{1,v,0}-Z_{s}^{2,v}\right\vert ^{2}\right) \mathrm{d}s\right] \\
&\leq &C\mathbb{E}^{\mathscr{F}_{t}}\left[ \int_{t}^{t+\delta }\left\vert
X_{s}^{1,v,0}-x\right\vert ^{2}\mathrm{d}s\right] \\
&\leq &C\delta \mathbb{E}^{\mathscr{F}_{t}}\left[ \sup_{t\leq s\leq t+\delta
}\left\vert X^{t,x;v,0}-x\right\vert ^{2}\right] \\
&\leq &\delta ^{2}.
\end{eqnarray*}%
Hence, it yields that%
\begin{eqnarray*}
\left\vert Y_{t}^{1,v,0}-Y_{t}^{2,v}\right\vert &=&\Bigg |\mathbb{E}^{%
\mathscr{F}_{t}}\Bigg [\int_{t}^{t+\delta }\Big (F\left(
s,X_{s}^{t,x;v,0},Y_{s}^{1,v,0},Z_{s}^{1,v,0},v_{s}\right) \\
&&-F\left( s,x,Y_{s}^{2,v},Z_{s}^{2,v},v_{s}\right) \Big )\mathrm{d}s\Bigg |
\\
&\leq &C\mathbb{E}^{\mathscr{F}_{t}}\left[ \int_{t}^{t+\delta }\left\vert
X_{s}^{1,v,0}-x\right\vert ^{2}\mathrm{d}s\right] \\
&&+C\delta ^{\frac{1}{2}}\mathbb{E}^{\mathscr{F}_{t}}\left[
\int_{t}^{t+\delta }\left\vert Y_{s}^{1,v,0}-Y_{s}^{2,v}\right\vert ^{2}%
\mathrm{d}s\right] ^{\frac{1}{2}} \\
&&+C\delta ^{\frac{1}{2}}\mathbb{E}^{\mathscr{F}_{t}}\left[
\int_{t}^{t+\delta }\left\vert Z_{s}^{1,v,0}-Z_{s}^{2,v}\right\vert ^{2}%
\mathrm{d}s\right] ^{\frac{1}{2}} \\
&\leq &C\delta ^{\frac{3}{2}}.
\end{eqnarray*}%
We thus complete the proof. ~\hfill $\Box $


\bigskip

\noindent \textbf{Acknowledgements. }{\small The author highly appreciates
the constructive comments and suggestions of the editor Prof. Alberto
Bressan and the referee which led to several improvements of the original
version. The author also wishes to thank Prof. Xiuqing Chen and Prof. Wulin
Suo for their valuable comments and discussions which improved the
presentation of this manuscript.}

\end{document}